\documentclass[final,leqno,onefignum,onetabnum]{siamltex1213}
\usepackage{amssymb}
\usepackage{amssymb}
\usepackage{amsmath,amssymb,graphicx}
\usepackage{multirow}
\usepackage{color}
\usepackage{epsfig}

\newcommand{\be}{\begin{eqnarray}}
\newcommand{\ee}{\end{eqnarray}}
\newcommand{\ben}{\begin{eqnarray*}}
\newcommand{\een}{\end{eqnarray*}}

\title{Behavior of different  numerical schemes for   population genetic drift problems}

\author{ Minxin Chen\thanks{Center for System Biology and Department of Mathematics,
Soochow University, Suzhou 215006, China (\email{chenmx@gmail.com}).}
\and  Chun Liu\thanks{Department of Mathematics, Pennsylvania State University, University Park, PA, 16802, USA (\email{liu@psu.edu}).}
\and Shixin Xu\thanks{Corresponding author. Department of Mathematics, Soochow
University, Suzhou 215006, China (\email{xsxztr@hotmail.com}). Department of Applied and Computational Mathematics and Statics, University of Notre Dame, Notre Dame, IN, 46556, USA.}
\and Xingye Yue\thanks{Department of Mathematics, Soochow
University, Suzhou 215006, China (\email{xyyue@suda.edu.cn}).}
\and Ran Zhang \thanks{Department of Mathematics, Fudan University, Shanghai 200433, China (\email{12110180047@fudan.edu.cn})}
 }

\begin{document}
\maketitle
\slugger{mms}{xxxx}{xx}{x}{x--x}%slugger should be set to mms, siap, sicomp, sicon, sidma, sima, simax, sinum, siopt, sisc, or sirev

\begin{abstract}
In this paper, we focus on numerical methods for the genetic drift problems, which is governed by a degenerated convection-dominated parabolic equation. Due to the degeneration and convection, Dirac singularities will always be developed at boundary points as time evolves. In order to find a \emph{complete solution} which should keep the conservation of total probability and  expectation, three different schemes based on  finite volume methods are used to solve the equation numerically: one is a upwind scheme, the other two are different central schemes.  We observed
that all the methods are stable and can keep the total probability, but have totally different long-time behaviors concerning with the conservation of expectation.   We prove that any extra infinitesimal diffusion leads to a same  artificial  steady state. So upwind scheme does not work due to its intrinsic numerical viscosity. We find one of the central schemes introduces a numerical viscosity term too, which is beyond the common understanding in the convection-diffusion community.  Careful analysis is presented to  prove that the other central scheme does work. Our study shows that for this kind of problems, the numerical methods must be carefully chosen and any method with intrinsic numerical viscosity should be avoided.
\end{abstract}



\pagestyle{myheadings}
\thispagestyle{plain}
\markboth{M. Chen, C. Liu, S. Xu, X. Yue and  R. Zhang}{ GENETIC DRIFT PROBLEM}

\section{Introduction}

The number of a particular gene (allele)  of one locus in the population varies randomly from generation to generation.
This process is a kind of stochastic process and named as genetic drift, which was  first introduced by one of the founders in the field of population genetics, S. Wright \cite{Wright1929}.
Genetic drift plays an important role in mutation, selection and molecular evolution \cite{Kimura1970,Kimura1983}.
Mathematical descriptions of genetic drift are typically built upon the Wright-Fisher model \cite{01,02} or its diffusion limit \cite{Moran1958,03, Kimura1962}.    The Wright-Fisher model describes the dynamics of   a gene with two alleles, A or B in a population of fixed size, N.  The model is formulated as a discrete-time Markov chain. If  $X_k$ denotes the number of A in $k_{th}$ generation,  then $N-X_k$ is the number of $B$ and the transition probabilities are the following binomial distribution:
\be
P(X_{n+1}=j|X_{n}=i)=C_N^j\left(\frac{i}{N}\right)^j\left(1-\frac{i}{N}\right)^{N-j}.
\ee
The first and second conditional moments of the Wright-Fisher process satisfy  \cite{Der}
\be\label{MFmoment}
E[X_{k+1}|X_{k}]=X_k,\\
Var[X_{k+1}|X_{k}]=X_k(1-\frac{X_k}{N}).
\ee
The first condition expresses the neutrality.
A basic phenomenon of neutral Wright-Fisher model without mutation is that, either all copies of allele A are lost from the population or   all the individuals carry allele A. This is the phenomenon of {\em fixation}.  And  the  allele will fix with probability equal to its initial frequency.    For a fixed population size, Wright-Fisher model works well, but if the size of the population changes over time, the diffusion approximation come to its own right.
In \cite{01} and  \cite{02}, the diffusion approximation is first introduced  to model the genetic drift, Moran \cite{Moran1958} and Kimura \cite{03} substantially extended and developed this approach. If we describe the process by $x(t)$, the fraction of type $A$ genes at time $t$, and $f(x,t)$ denotes the probability density of $x$ at time $t$,
Kimura  \cite{03, Kimura,zhao} showed that  $f(x,t)$ obeys the following diffusion equation,
\begin{equation}
\frac{\partial f(x,t)}{\partial t}+\frac{\partial j(x,t)}{\partial
x}=0, \nonumber
\end{equation}
where the quantity $j(x,t)$ is the current that characterizes the flow of probability density, with the form as
\begin{equation}
j(x,t)=R(x)f(x,t)-\frac{1}{4{\it N}}\ \frac{\partial}{\partial x}(x(1-x)f(x,t)). \nonumber
\end{equation}
$R(x)$ represents the deterministic part of gene frequency dynamics and is typically taken as a polynomial in {\it x} whose coefficients depend on mutation rates,\ migration rates and selection coefficients.
At the first step,\ we will consider the simplest case,
 where the only evolutionary force acting on a randomly mating diploid population is diffusion,
thus $R(x)=0$.
 By rescaling of the time, we can get the model introduced by Crow and Kimura \cite{Kimura1970},
\begin{equation}
\frac{\partial f(x,t)}{\partial t}+\frac{\partial j(x,t)}{\partial x}=0,\label{con-dif0}
%\frac{\partial f(x,t)}{\partial t}-\frac{\partial}{\partial x}\left(x(1-x)\frac{\partial f}{\partial x}\right)+\frac{\partial}{\partial x}((2x-1)f)=0, \ \ x\in (0,1),
\end{equation}
where  the quantity $j(x,t)$ is the current that characterizes the flow of probability density, with the form as
\begin{equation}
j(x,t)=-\frac{\partial}{\partial x}\left(x(1-x)f(x,t)\right),\label{eq03}
\end{equation}
and the  zero   current boundary conditions   as (\cite{05})
\begin{eqnarray}
&&j(0,t) = 0,\ \ \ \ j(1,t) = 0. \label{bc1}
\end{eqnarray}

% \begin{remark}
%Here we only consider the simple case and the mutations are not included. If we consider the mutation, such as in \cite{05} and \cite{zhao}, the flux is given in the following form
%\begin{equation}
%j(x,t)=R(x)f(x,t)-  \frac{\partial}{\partial x}(x(1-x)f(x,t)),\label{current}
%\end{equation}
%where $R(x)$ represents the deterministic part of gene frequency dynamics and is typically taken as a polynomial in {\it x} whose coefficients depend on mutation rates,\ migration rates and selection coefficients.
%\end{remark}

 %The previous solutions to the equation,\ such as the level of
%heterozygosity in the population (see \cite{04}),\ suffer from a
%lack of completeness.\ It is taken to hold only on the interior of
%the possible range of {\it x} but not at the boundary values of
%{\it x}, i.e.,\ not at $x=0$ and $x=1$.

In \cite{Tran2013}, Tran et al. defined that a distrbution $f\in H$ is called as  a weak solution of problem \eqref{con-dif0}-\eqref{bc1} with a initial state $f(x,0)= f_0(x) \in H$, if  \be\label{weak}
(f_t,\phi)=(f,x(1-x)\frac{\partial^2 \phi}{\partial x^2}), ~~ \forall \phi\in C^{\infty}[0,1],
\ee
where  $H=\{f:[0,1]\rightarrow [0,\infty]| \int_0^1fg dx<\infty,~ \forall g\in C^{\infty}[0,1]\}$ is the set of all general distribution functions on $[0,1]$.
They also proved   the existence and uniqueness of the weak solution.

 In the definition \eqref{weak}, if we set   $\phi=1$, we can find that  problem \eqref{con-dif0}-\eqref{bc1} conserve the total probability, for any time $t>0$,
\begin{equation}
 \int_{0}^{1}\ f(x,t)\ \mathrm d x =  \int_{0}^{1}\ f(x,0)\ \mathrm d x = 1\label{conp}
\end{equation}
Furthermore, if we chose   $\phi=x$ in \eqref{weak}, then we get the conservation of   the mean gene frequency (expectation), i.e.   \be
\int_0^1 xf(x,t)\
\mathrm{d}x\equiv \int_0^1 xf(x,0)\ \mathrm{d}x.\label{expp}\ee

McKane and Waxman (\cite{05}) gave a closed form of the
singular solution to the above genetic drift problem as
\begin{equation}
f(x,t)={\Pi}_0(t) \delta(x)+{\Pi}_1(t) \delta(1-x)+f_r(x,t),\label{eq}
\end{equation}
where $\Pi_0(t),~\Pi_1(t)$ and $f_r(x,t)$ are smooth functions and $f_r(x,t)$ is in a form of infinite series. \eqref{eq} means the solution has three parts: two singular part ${\Pi}_0(t) \delta(x),~{\Pi}_1(t) \delta(1-x)$ and one smooth part $f_r$.
Specifically,  if we solve this equation \eqref{con-dif0} subject to  the condition that all replicate populations initially have the gene frequency of $p$, so $f(x,0)$ corresponds to an initial distribution where only the single frequency $p$ is present,    \begin{equation} \label{ic} f(x,0)=\delta(x-p), \ \ 0<p<1,\end{equation}
which is a Dirac delta function at $x=p$.  In this case, by the conservation of total probability \eqref{conp} and the conservation of expectation \eqref{expp},   Mckane and Waxman proved that
\begin{equation}\lim_{t\rightarrow\infty}f(x,t)=(1-p)\delta(x)+p\delta(1-x)\label{stable}\end{equation}
in the sence of distribution. This is just the  {\em fixation phenomenon} corresponding to the Wright-Fisher model.
However, it is not easy to compute the infinite
series  in (\ref{eq}), and for more complex situations, e.g. $R(x)\not\equiv 0$, one can hardly derive the  closed form solutions as above, thus direct numerical simulation is also necessary for this problem.

 The numerical scheme should be design  to find a \emph{complete solution} \cite{zhao}, which should keep the conservations of total probability \eqref{conp} and  expectation \eqref{expp}.
Among various numerical methods, we choose finite volume method (FVM)  \cite{Chen2013,Eymard,Roos1996} since it is easy to keep the conservation laws numerically.
 \eqref{con-dif0} looks like a pure diffusion equation. Actually it is a convection-dominated diffusion equation. This can be seen if we  rewrite \eqref{con-dif0} as
\begin{equation}
\frac{\partial f(x,t)}{\partial t}-\frac{\partial}{\partial x}\left(x(1-x)\frac{\partial f}{\partial x}\right)+\frac{\partial}{\partial x}((2x-1)f)=0  \label{origin}
\end{equation}
with diffusion coefficient
$ x(1-x)$ and convection velocity
$ (2x-1)$. It is clear that (\ref{origin})
is convection-dominated near boundaries $x=0$ and $x=1$, where convection
velocity is up to 1 while diffusion coefficient is down to $0$.

%\begin{remark}
 At the boundary points,  \eqref{origin}  is  degenerated  to a pure convection problem.   Wether or not a boundary condition should be imposed is the key point to obtain a well-posed solution.  For problem \eqref{origin}, the situation is that we can get a unique regular solution without any boundary condition  imposed. But this solution is not a complete one: it does not keep the conservation of probability \eqref{conp}. So the no-flux boundary condition \eqref{bc1} is imposed and of course we can not expect to obtain a regular solution now. See the discussions  on the   boundary condition  in the appendix of \cite{zhao}.
%\end{remark}

 In general, a
upwind FVM (UFVM) scheme is a better choice for the convection-dominated diffusion
problem to achieve stability due to its intrinsic numerical
viscosity. It does give us a stable  numerical solution. But we find that we always get the same long-time behavior no matter what the initial state is.  This is obviously wrong for gene drift (See \eqref{stable}) because it can not keep the conservation of the expectation \eqref{expp}.

To see what's wrong here, we appeal to  the method of vanishing viscosity,\ i.e.,\ a
 small viscosity term is first added, then
the limit behavior of the solution is considered when  the added viscosity
tends to zero. \ We see that the limitation of the steady state solution  is uniquely determined and thus it has nothing to do with the initial conditions. This means that the long time behavior of the original
problem will be changed by any added infinitesimal viscosity.\ That is
the reason why the upwind scheme does not work
 for the genetic drift problem.

When we turn to the central difference  FVM, we have two choices: in the first method,  the central difference is applied to   discrete the fluxes  induced by  the diffusion $-\frac{\partial}{\partial x} (x(1-x)\frac{\partial f}{\partial x} )$ and convection $\frac{\partial}{\partial x}((2x-1)f)$ respectively (this is referred to Scheme 2 in next section); in the other scheme,  we keep the flux \eqref{eq03} as a whole, which is also discretized  by the central difference (this is referred to Scheme 3).  We find that both Schemes 2 and 3 are stable, which is a surprising since we are solving a convection-dominated problem by central schemes. We also observe that the steady-state of Scheme 2 is the same as UFVM (ref. to Scheme 1 in next section) but it takes a much longer time for Scheme 2 to achieve the steady state  than for Scheme 1.
Scheme 3 is the simplest one and we find it gives us a \emph{complete  solution}. Dirac singularities  develop at both boundary points with proper weights rather than for Schemes 1 and 2, Dirac singularities with same weight always develop at both ends.  For an explanation, we give a careful analysis to show that Scheme 3  always converges to the right solution and Scheme 2 is equivalent to Scheme 3 plus a 2nd order $O(h^2)$ viscosity term though it is a central difference scheme itself. That is the reason why it takes a much longer time for Scheme 2 to achieve the steady state  than for Scheme 1, which, as we all have known,  is equivalent to Scheme 2 plus a much larger first order $O(h)$ viscosity term.

 %Two finite volume methods, upwind FVM (UFVM) and central FVM (CFVM)
%schemes are used to solve the equation numerically. We observed
%that the long time behaviors of the numerical solutions of these two
%methods are totally different. Based on the conservations of total
%probability and the mean gene frequency (expectation), the conclusion is drawn
%that the results of Scheme 3 are correct while the
%results of   UFVM are not correct since it destroys
%the conservation of the mean gene frequency. However,

 This paper is organized as follows. We present three different finite volume schemes for
 genetic drift problem in Section 2. Numerical results and
 analysis for three methods are presented in Section 3. The final section  gives some concluding remarks and discussions.

\section{Numerical Methods}
\setcounter{equation}{0}
%By re-scaling in time, we can always replace   the coefficients $\frac{1}{4N}$ in \eqref{eq03origin} and \eqref{origin} by $1$, i.e.
%\begin{equation}
%\frac{\partial f(x,t)}{\partial t}+\frac{\partial j(x,t)}{\partial x}=0,\label{con-dif0}
%%\frac{\partial f(x,t)}{\partial t}-\frac{\partial}{\partial x}\left(x(1-x)\frac{\partial f}{\partial x}\right)+\frac{\partial}{\partial x}((2x-1)f)=0.
%\end{equation}
%where  \begin{equation}
%j(x,t)=-\frac{\partial}{\partial x}\left(x(1-x)f(x,t)\right)=-\left((x(1-x))\frac{\partial f}{\partial x}+(1-2x)f\right).\label{eq03}
%\end{equation}
In order to keep the total probability,    we start from  FVM. A uniform grid,\ with grid spacing $h=1/M$
and grid points $x_k=ih,\ i=0,\dots,M,$ is used to discretize the
space domain $[0,1]$. Likewise,\ the time domain is uniformly
discretized   with  step size $\tau$.\ Let $j_i^n$ and $f_i^n$ be the numerical approximations of
$j(x_i,t_n)$,
$f(x_i,t_n)$, respectively. For inner mesh point $x_i$ $(1\le i\le M-1)$, the
control volume is
$$
\mathcal{D}_{i}=\{x|\ x_{i-\frac{1}{2}}\le x \le
x_{i+\frac{1}{2}}\},
$$
 where $x_{i+\frac{1}{2}}=(i+\frac{1}{2})h$.\
By FVM, we have
\begin{equation}
\frac{f_ i^{n+1}-f_i^n}{\tau}+\frac{j_{i+\frac{1}{2}}^{n+1}-j_{i-\frac{1}{2}}^{n+1}}{h} = 0.\label{eq8}
\end{equation}

For the   boundary points $x_0=0$ and $x_M=1$, the control volumes are
$$
\mathcal{D}_{0}=\{x|\ 0\le x \le
x_{\frac{1}{2}}\}, \ \mbox{and}\ \mathcal{D}_{M}=\{x|\ x_{M-\frac{1}{2}}\le x \le
1\}.
$$
By the boundary conditions $j(0,t)=j(1,t)=0$,
%and
%\begin{equation}
%\int_{x_0}^{x_\frac 1 2}\ f(x,t_{n})\  \mathrm{d} x \approx \frac{h}{2}f_{0,n},
%\end{equation}
%\begin{equation}
%\int_{t_n}^{t_{n+1}}\ j(x_{\frac{1}{2}},t) \ \mathrm{d} t \approx \tau j(x_\frac{1}{2},t_{n+1}),
%\end{equation}
%\begin{equation}j( x_{\frac{1}{2}},t_{n+1})=(2x_{\frac{1}{2}}-1)f(x_{\frac{1}{2}},t_{n+1})-
%x_{\frac{1}{2}}(1-x_{\frac{1}{2}})
%  \frac{\partial f}{\partial x}\Big{|}_{\frac{1}{2},n+1},
% \end{equation}
%\begin{equation} \frac{\partial f}{\partial x}\Big{|}_{\frac{1}{2},n+1}  \approx
%\frac{f_{0,n+1}-f_{1,n+1}}{h},
%\end{equation}
we have
\begin{equation}\label{eq6}
\frac{f_ 0^{n+1}-f_0^n}{\tau}+\frac{j_{\frac{1}{2}}^{n+1}-0}{h/2} = 0 \ \mbox{and}\
\frac{f_M^{n+1}-f_M^n}{\tau}+\frac{0-j_{M-\frac{1}{2}}^{n+1}}{h/2} = 0.
\end{equation}
%\begin{equation}
%\frac{h}{2}(f_{0,n+1}-f_{0,n})=\tau\left((2x_{\frac{1}{2}}-1)f(x_{\frac{1}{2}},t_{n+1})-
%x_{\frac{1}{2}}(1-x_{\frac{1}{2}})\frac{f_{0,n+1}-f_{1,n+1}}{h}\right).\label{eq6}\end{equation}
%Let $j=0$ in Eq.(9) and (11),\ we get $ j(x_{\frac 1 2},t_{n+1})$.\ Substituting Eq.(16)-- (18) into Eq.(15),\ we get
%\begin{equation}
%f(x_0,t_{n+1})-f(x_0,t_n)=-2\nu[b_0f(x_{0},t_{n+1})+c_0f(x_1,t_{n+1})]
%\end{equation}
%where
%\begin{eqnarray}
%b_0&=\frac 1 2(2x_{\frac 1 2}-1)+\frac 1 h x_{\frac {1}{2}}(1-x_{\frac {1}{2}})
%\\
%c_0&=\frac 1 2(2x_{\frac 1 2}-1)-\frac 1 h x_{\frac {1}{2}}(1-x_{\frac {1}{2}})
%\end{eqnarray}
%At the right boundary point $x=1$, \eqref{eq1} is discretized analogously,
%\begin{equation}
%\frac{h}{2}(f_{M,n+1}-f_{M,n})=\tau
%j_{M-\frac{1}{2},n+1}\label{eq9}
%\end{equation}

To get a  fully discrete  scheme, we still need to approximate the
term $j_{i+\frac{1}{2}}^{n+1}$, for $i=0,...,M-1$. It is
treated differently by the following three FVM schemes:
\begin{itemize}
 \item Scheme 1:  approximate $j(x ,t)=- x (1-x )\frac{\partial f}{\partial x}+(2x-1)f$ at point $x_{i+\frac 1 2}$
  by   upwind scheme
\begin{eqnarray}\ \ \displaystyle
\begin{footnotesize}{j_{i+\frac 1 2}^{n+1} \!\!=\!\! \left\{\begin{array}{l}\!\!\! -x_{i+\frac 1 2}(1-  x_{i+\frac 1 2})\frac{f_{i+1}^{n+1}-f_{i}^{n+1}}{h}+ (2x_{i+\frac 1 2 }-1)f_{i+1}^{n+1},~ 2x_{i+\frac 1 2 }-1<0,\\
 \!\!\! -x_{i+\frac 1 2}(1-x_{i+\frac 1 2})\frac{f_{i+1}^{n+1}-f_{i}^{n+1}}{h} +(2x_{i+\frac 1 2 }-1)f_{i}^{n+1},~  2x_{i+\frac 1 2 }-1>0.\end{array}\right. \label{scheme0}}\end{footnotesize}
\end{eqnarray}
   \item Scheme 2:   approximate $j(x ,t)=- (x (1-x ))\frac{\partial f}{\partial x}+(2x-1)f $ at point $x_{i+\frac 1 2}$
  by   central scheme
    \begin{equation}
j_{i+\frac{1}{2}}^{n+1} =-
 x_{i+\frac{1}{2}}(1-x_{i+\frac{1}{2}})\frac{f_{i+1}^{n+1}-f_{i}^{n+1}}{h}+( 2x_{i+\frac 1 2}-1)\frac{f_{i+1}^{n+1}+f_i^{n+1}} 2.\label{scheme1}
\end{equation}

 \item Scheme 3: approximate $j(x,t)=-\frac{\partial}{\partial x}\left(x(1-x)f(x,t)\right)$ at point $x_{i+\frac 1 2}$ by  central scheme
  \begin{equation}
j_{i+\frac{1}{2}}^{n+1} =-\frac{ x_{i+1}(1-x_{i+1})f_{i+1}^{n+1}-x_{i}(1-x_{i})f_{i}^{n+1}}{h}.\label{eq7}
\end{equation}
 \end{itemize}

Scheme $3$ was recently used in \cite{zhao} for a numerical investigation on the random genetic problems, where some applications can be found on more complicated  topics such as time-dependent probability of fixation for a neutral locus or in the presence of selection effect within a population of constant size; probability of fixation in the presence of selection and demographic change. In this paper, we confine ourself to the simplest case to see the behaviors of different schemes.

\section{Numerical Results and Analysis}
\setcounter{equation}{0}
\subsection{Numerical Results}
\begin{figure}[!htbp]
\begin{center}
\includegraphics[width=2.in]{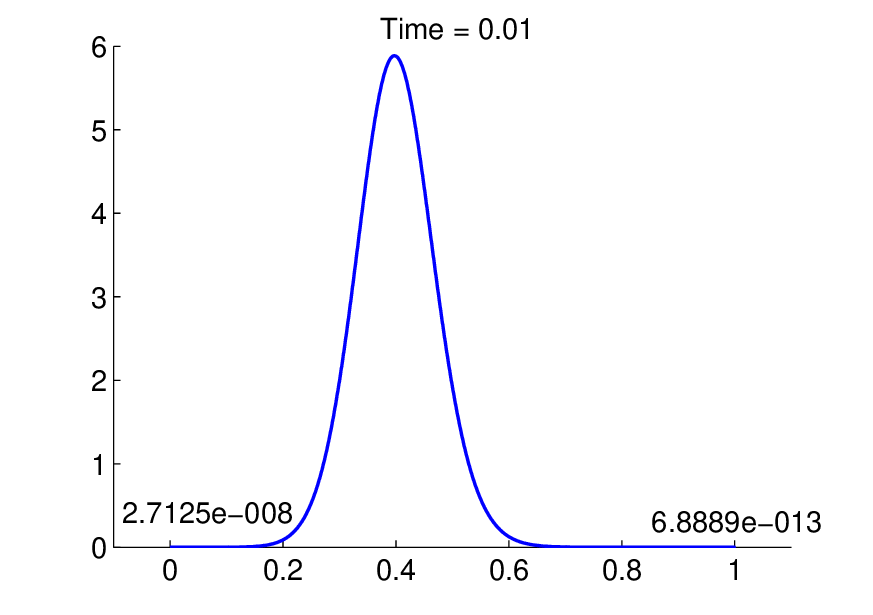}
\includegraphics[width=2.in]{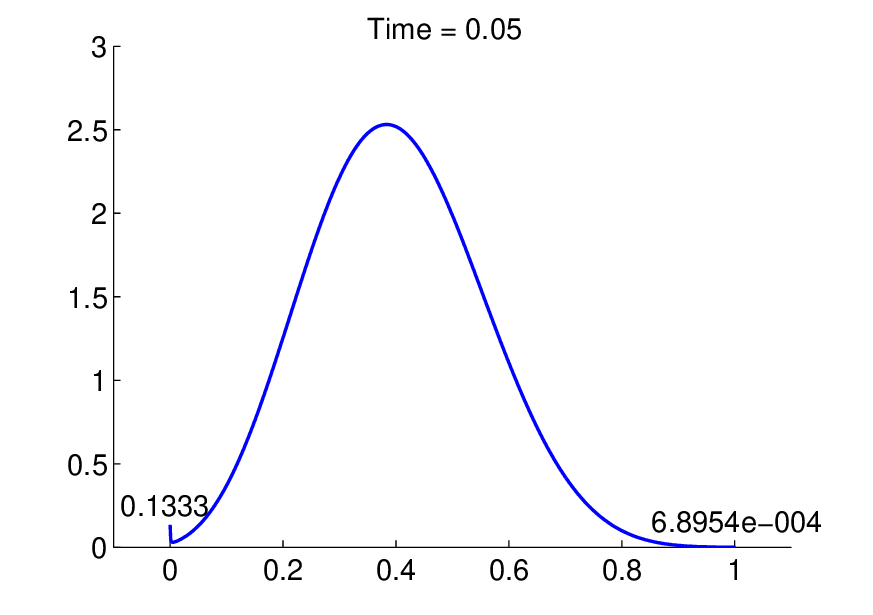}\\
\includegraphics[width=2.in]{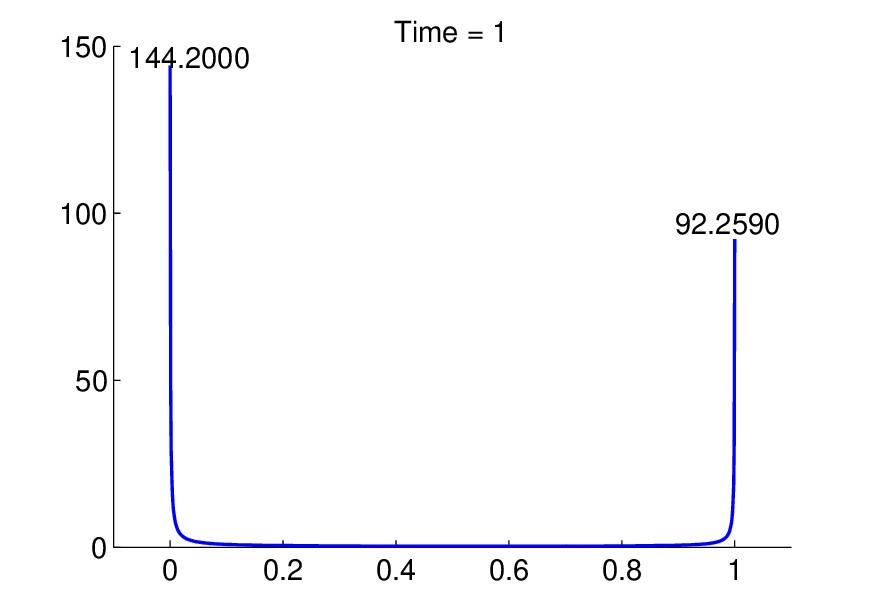}
\includegraphics[width=2.in]{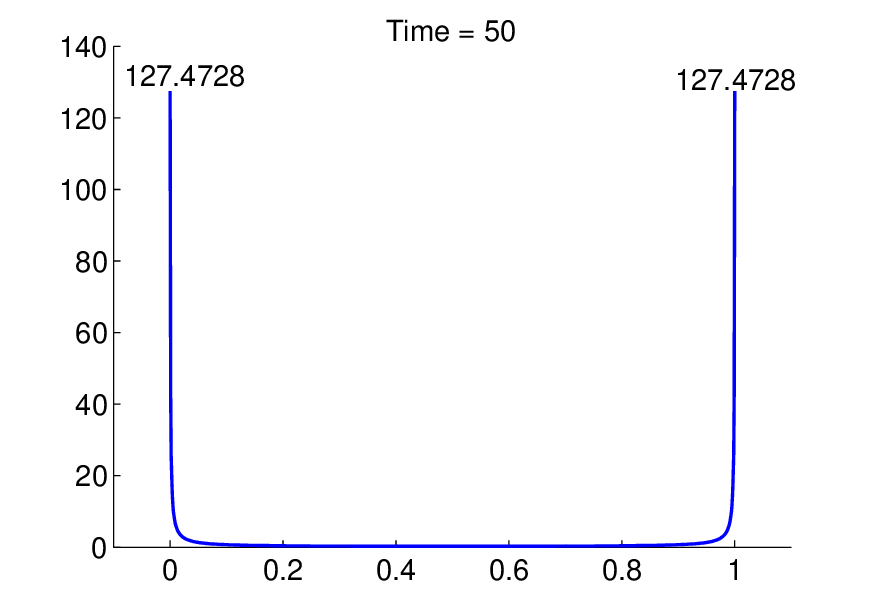}
\end{center}
\caption{Numerical results of Scheme 1  with initial state $f(x,0) \sim
N(0.4,0.01^2)$ at different time  $t= 0.01, 0.05, 1, 50$. The   step sizes are $h=1/1000,\ \tau =1/1000$. }\label{up04}
\end{figure}

\begin{figure}[!htbp]
\begin{center}
\includegraphics[width=2.in]{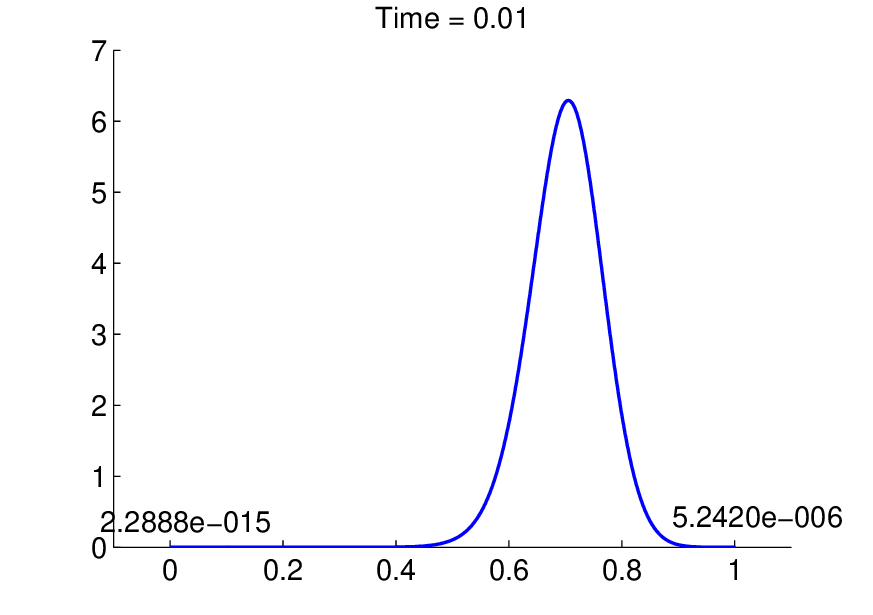}
\includegraphics[width=2.in]{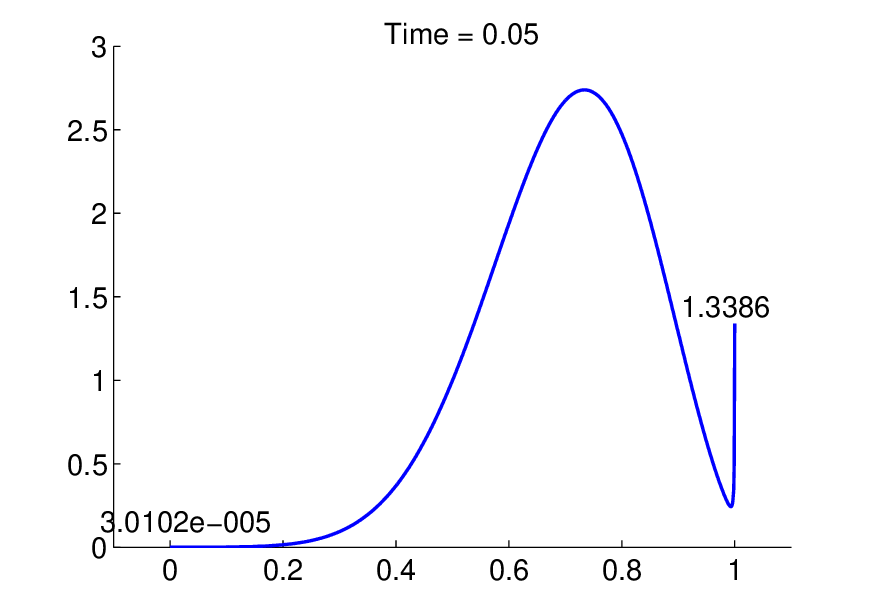}\\
\includegraphics[width=2.in]{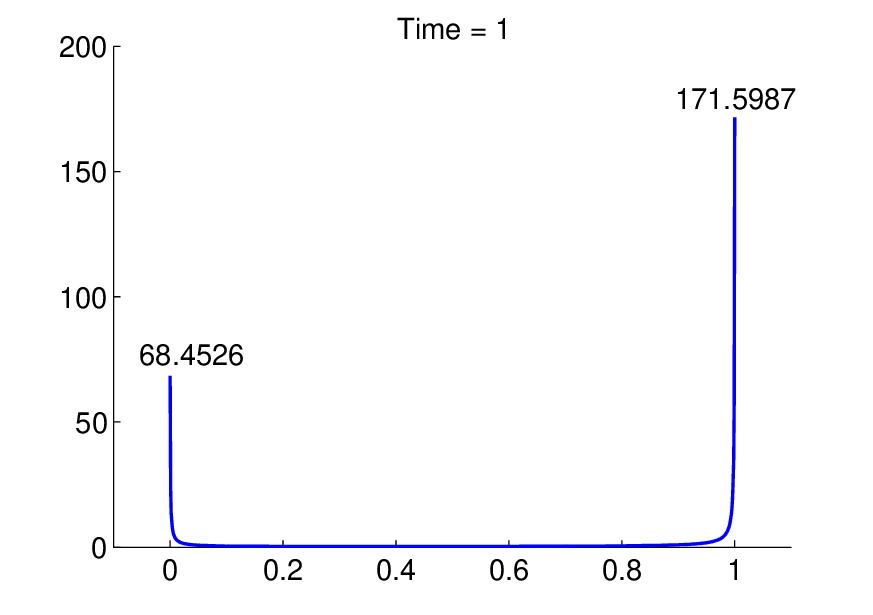}
\includegraphics[width=2.in]{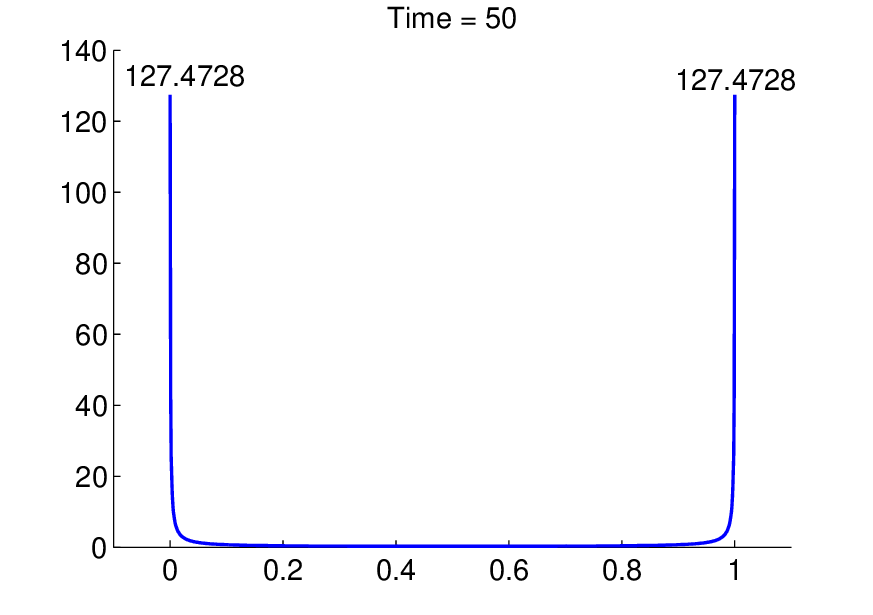}
\end{center}
\caption{Numerical results of Scheme 1  with initial state $f(x,0) \sim
N(0.7,0.01^2)$ at different time $t= 0.01, 0.05, 1, 50$. The step sizes are $h=1/1000,\ \tau =1/1000$. }\label{up07}
\end{figure}

\begin{figure}[!htbp]
\begin{center}
\includegraphics[width=2.in]{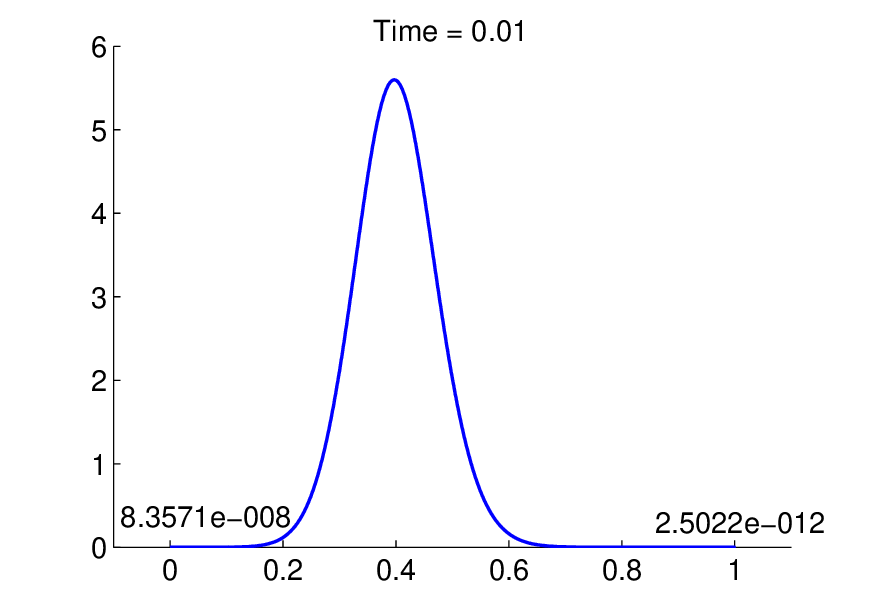}
\includegraphics[width=2.in]{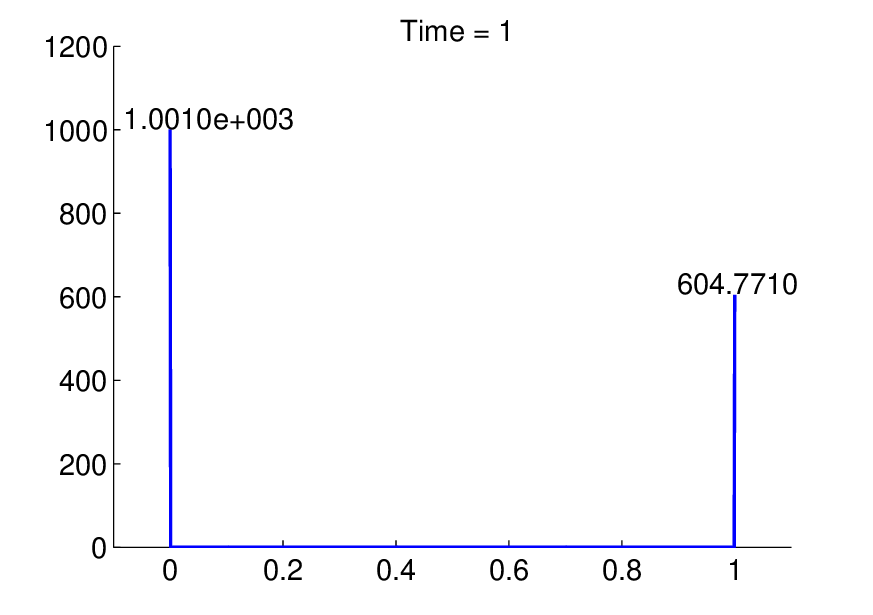}\\
\includegraphics[width=2.in]{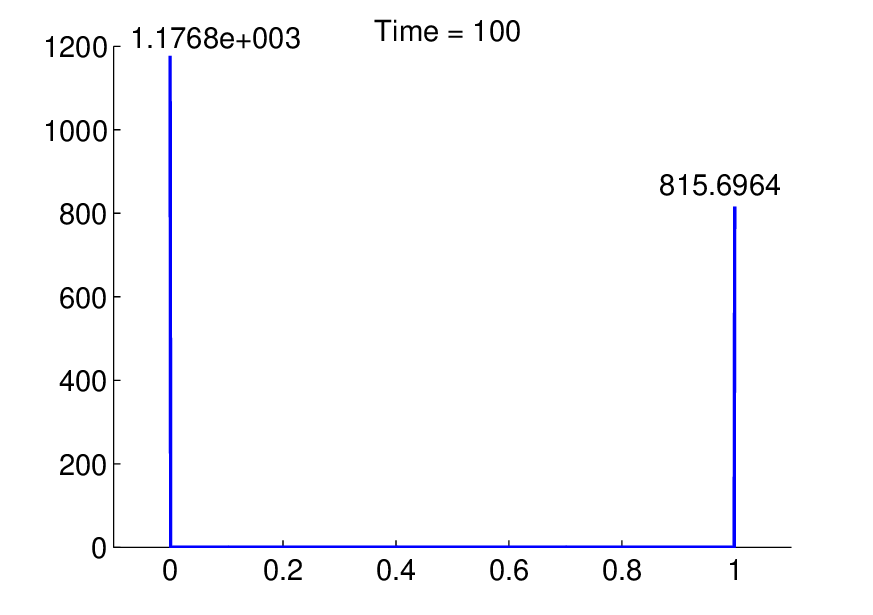}
\includegraphics[width=2.in]{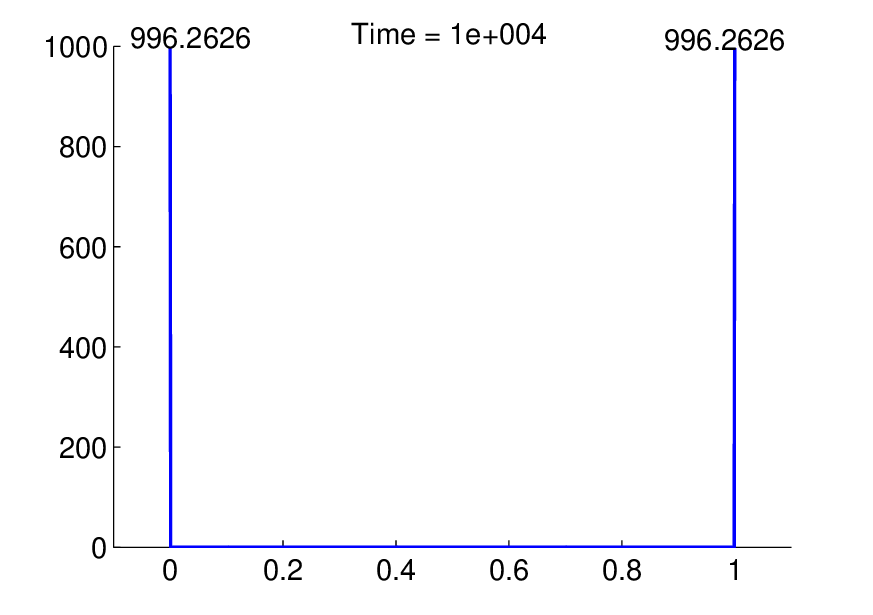}
\end{center}
\caption{Numerical results of Scheme 2 with initial state $f(x,0) \sim
N(0.4,0.01^2)$ at different time $t= 0.01, 1, 100, 10000$. The step sizes are $h=1/1000,\ \tau =1/1000$. }\label{s04}
\end{figure}

\begin{figure}[!htbp]
\begin{center}
\includegraphics[width=2.in]{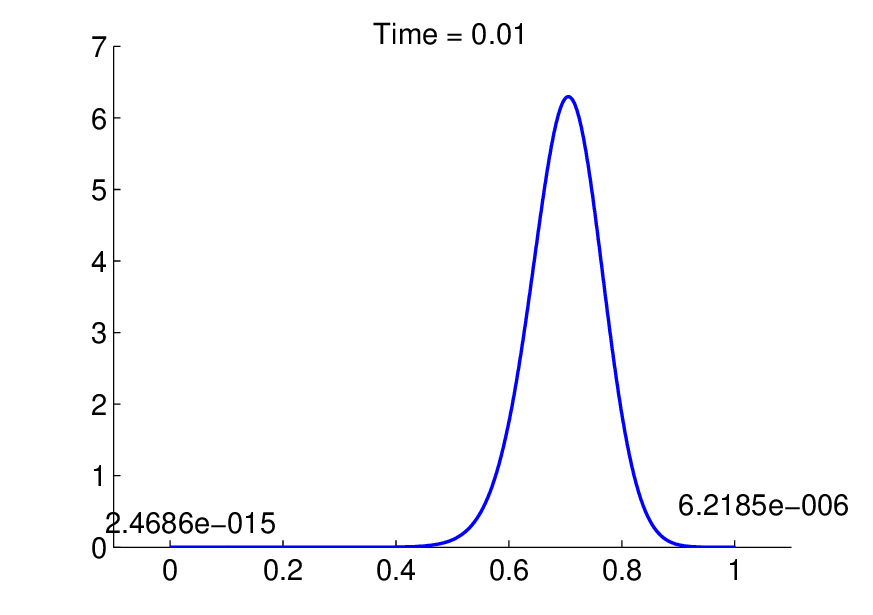}
\includegraphics[width=2.in]{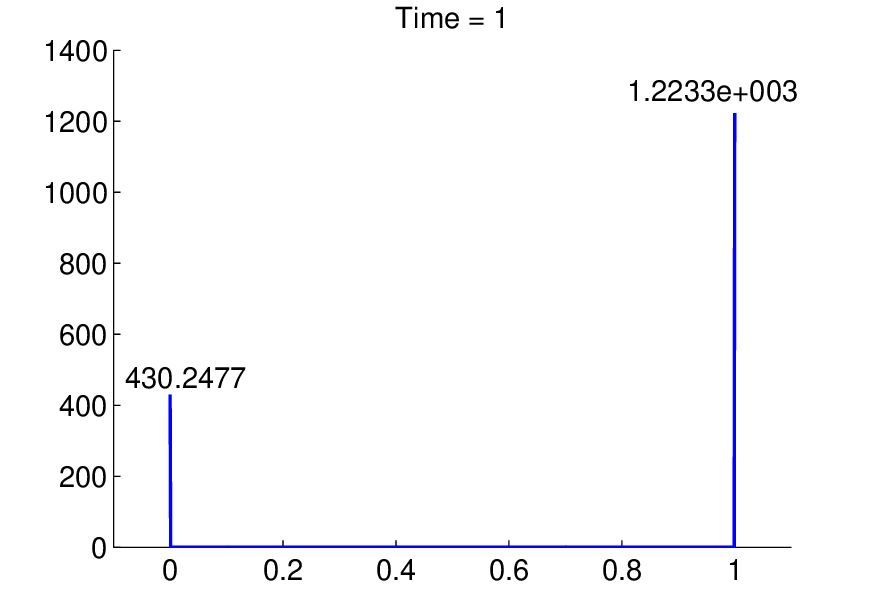}\\
\includegraphics[width=2.in]{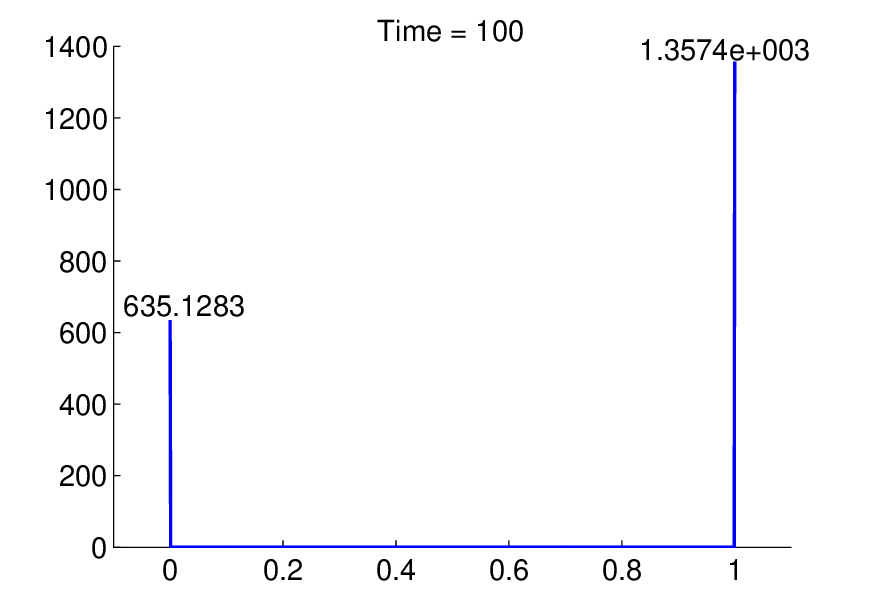}
\includegraphics[width=2.in]{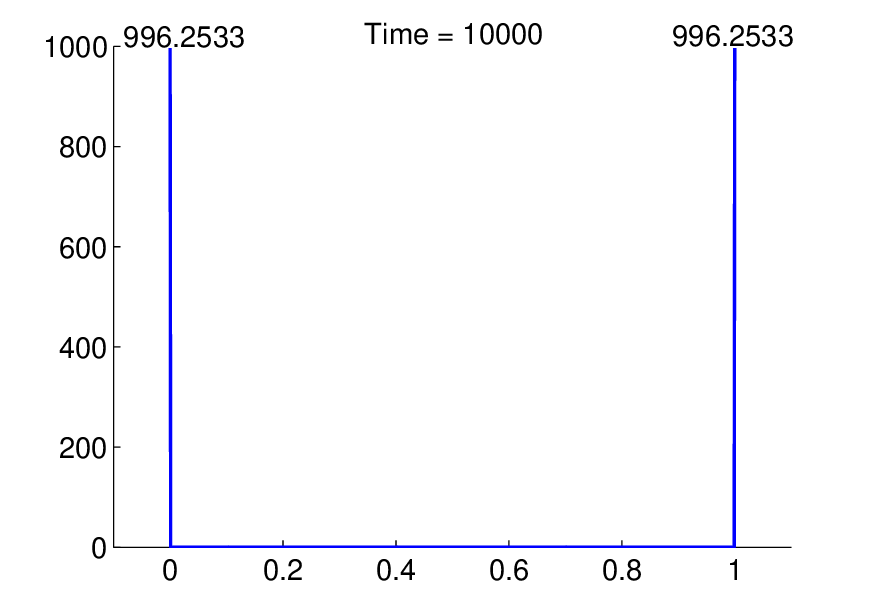}
\end{center}
\caption{Numerical results of Scheme 2 with initial state $f(x,0) \sim
N(0.7,0.01^2)$ at different time $t= 0.01, 1, 100, 10000$. The step sizes are $h=1/1000,\ \tau =1/1000$. }\label{s07}
\end{figure}

\begin{figure}[!htbp]
\begin{center}
\includegraphics[width=2.in]{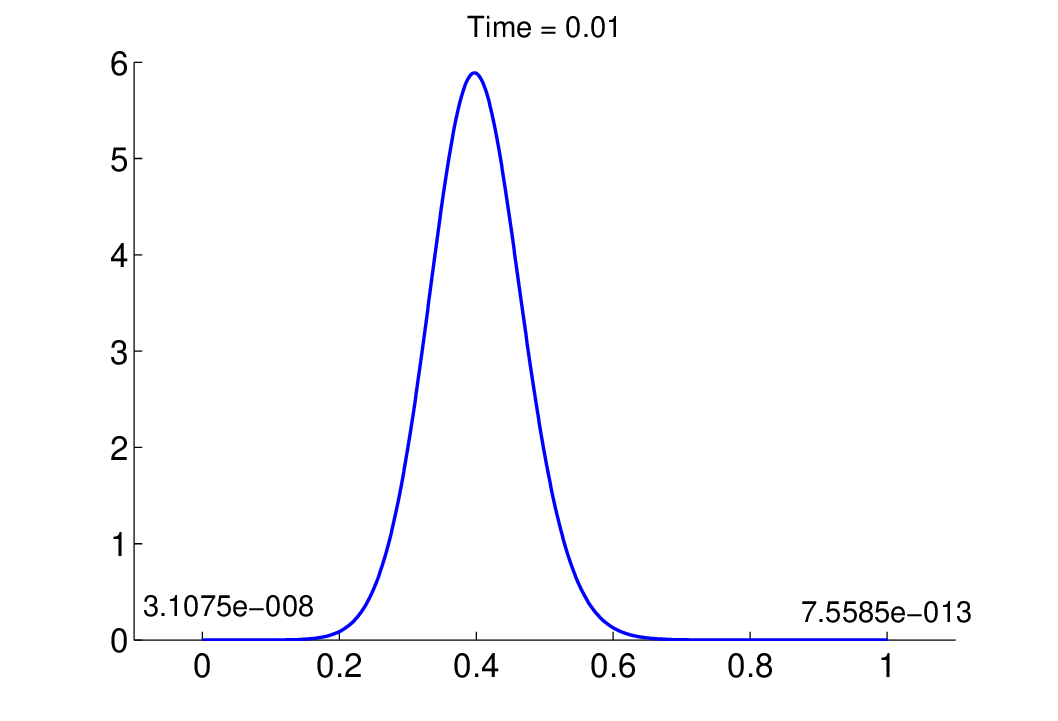}
\includegraphics[width=2.in]{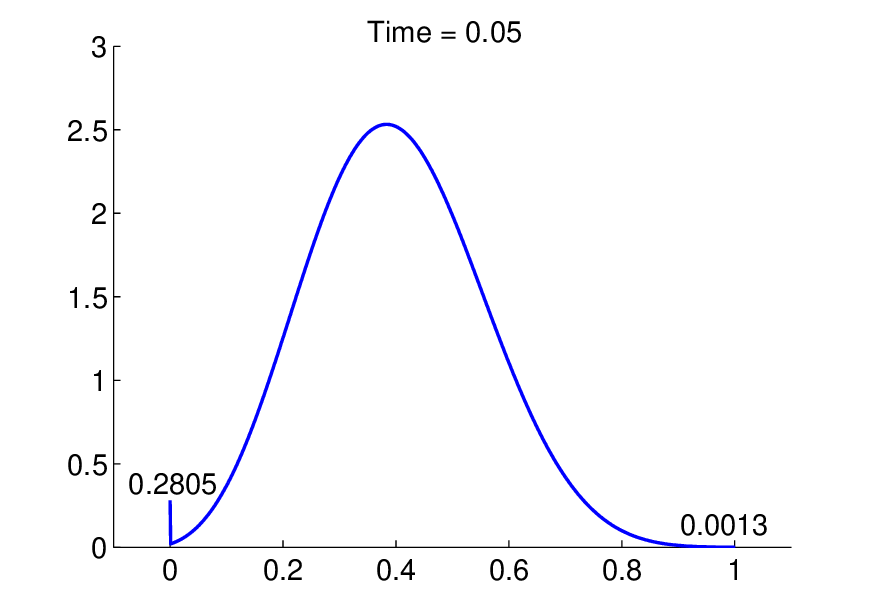}\\
\includegraphics[width=2.in]{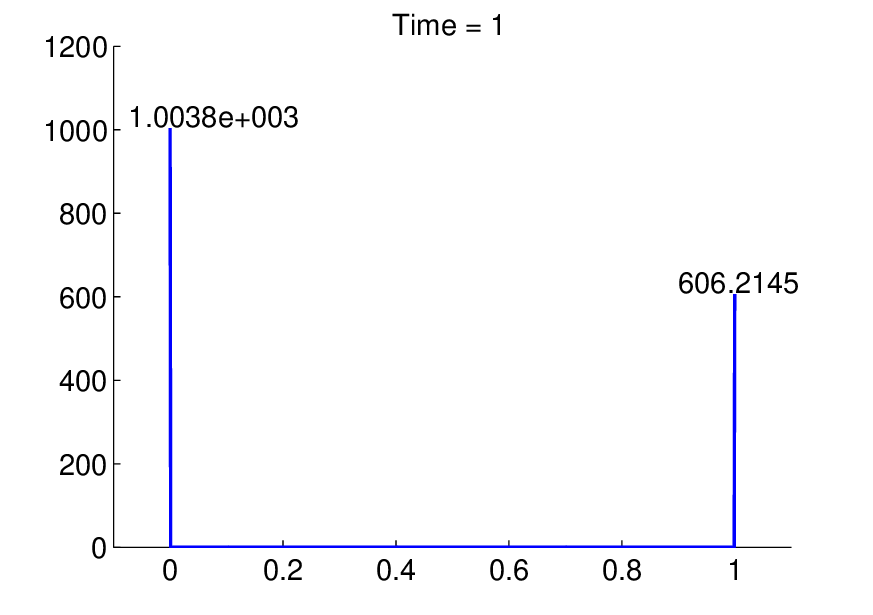}
\includegraphics[width=2.in]{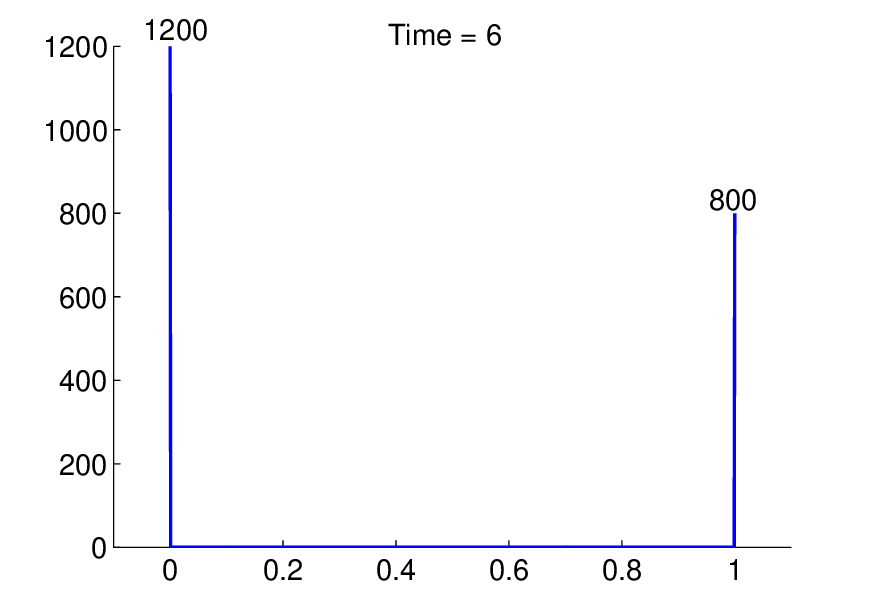}\\
\end{center}
\caption{Numerical results of Scheme 3 with initial state $f(x,0) \sim
N(0.4,0.01^2)$ at different time $t=0.01, 0.05, 1, 6$. The step sizes are $h=1/1000,\ \tau =1/1000$. }\label{c04}
\end{figure}

\begin{figure}[!htbp]
\begin{center}
\includegraphics[width=2.in]{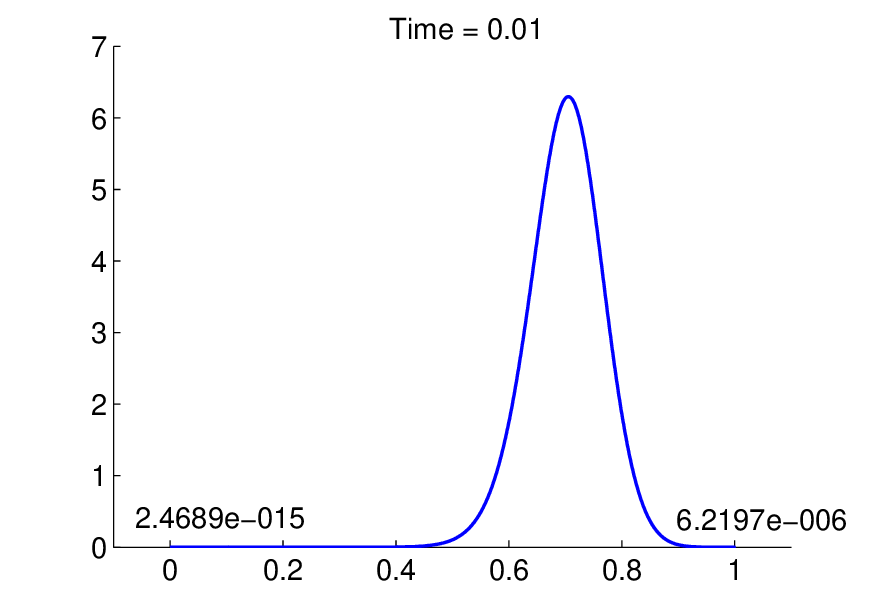}
\includegraphics[width=2.in]{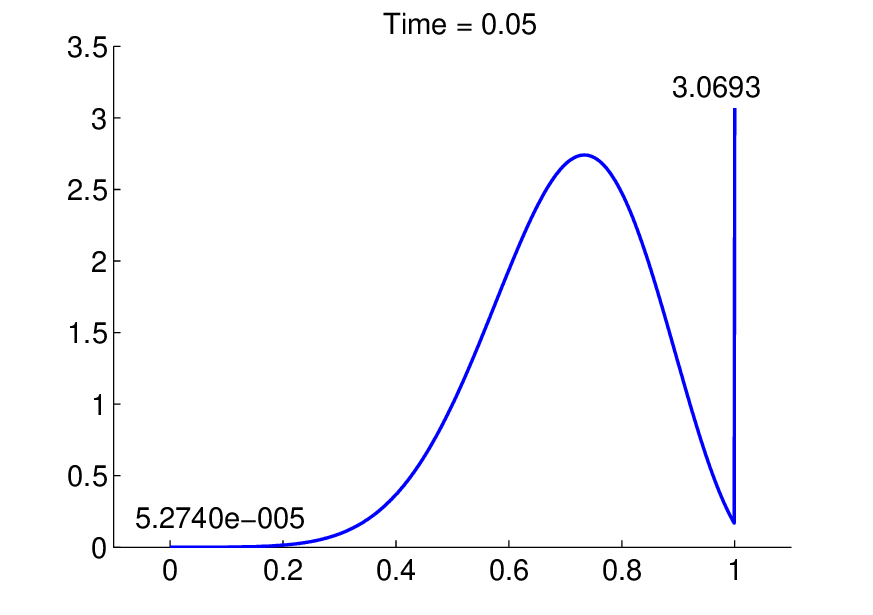}\\
\includegraphics[width=2.in]{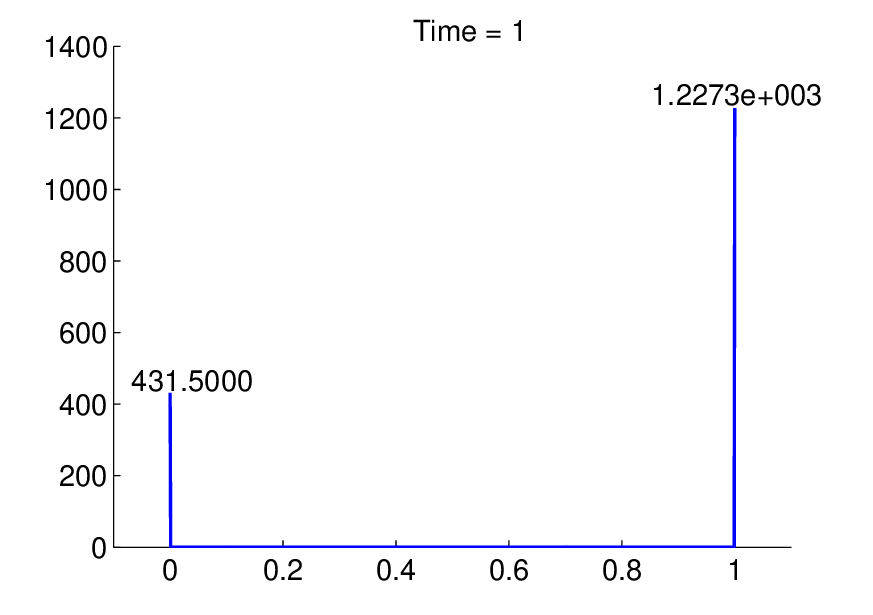}
\includegraphics[width=2.in]{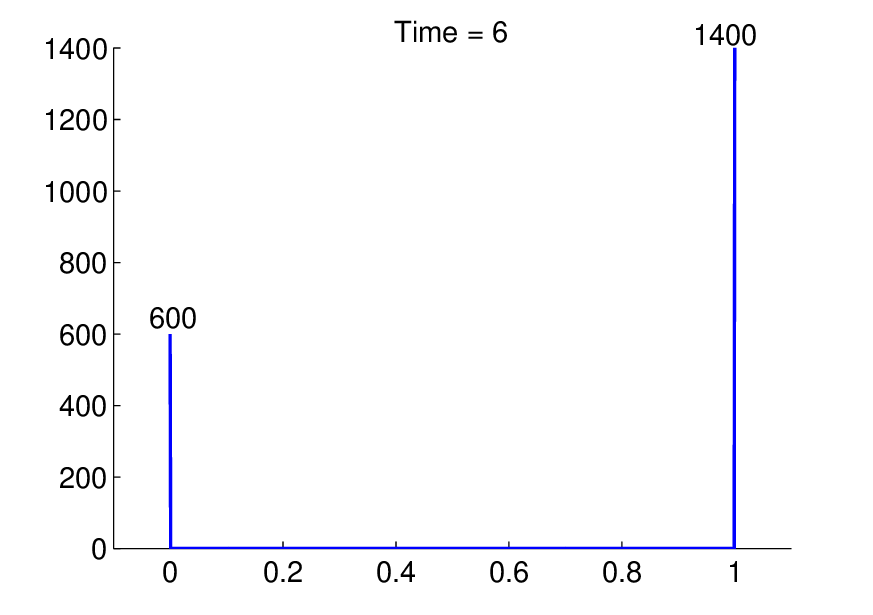}\\
\end{center}
\caption{Numerical results of Scheme 3 with initial state $f(x,0) \sim
N(0.7,0.01^2)$ at different time $t=0.01, 0.05, 1, 6$. The  step sizes are $h=1/1000,\ \tau =1/1000$. }\label{c07}
\end{figure}

The numerical results of different schemes  for the
genetic drift problem  \eqref{con-dif0}, \eqref{eq03}, \eqref{bc1}  and \eqref{ic}   are shown in this section.   In the numerical simulation,
we first approximate the initial state Dirac delta function \eqref{ic} by a
normal distribution function:
\begin{equation}
f(x,0)\sim N(p,{\sigma}^2)
\end{equation}
with $\sigma << 1.$ Later,\ we will show the dependence on the
mean value $p$ of the results.
In Figs. \ref{up04}-\ref{c07} the numerical probability density at
different time with various space grid sizes and initial
states, are shown for all three methods. These figures show that, for all schemes,
the probability density vanishes in the interval (0,1) and forms two peaks at the boundary ($x=0$
and $x=1$) as time evolves. But the heights of the peaks at the
boundary at the steady state are different. For Scheme 3, the heights of the two peaks at
the boundary of steady-state solution are depended on the mean of
initial probability density. For Schemes 1 and 2, the heights are
equals and has no relations to the expectation of initial condition. This means Schemes 1 and 2 can not keep the conservation of expectation and  can not give a \emph{complete solution}.

\begin{figure}[!htbp]
\begin{center}
\includegraphics[height=2.1in, width=2.5in]{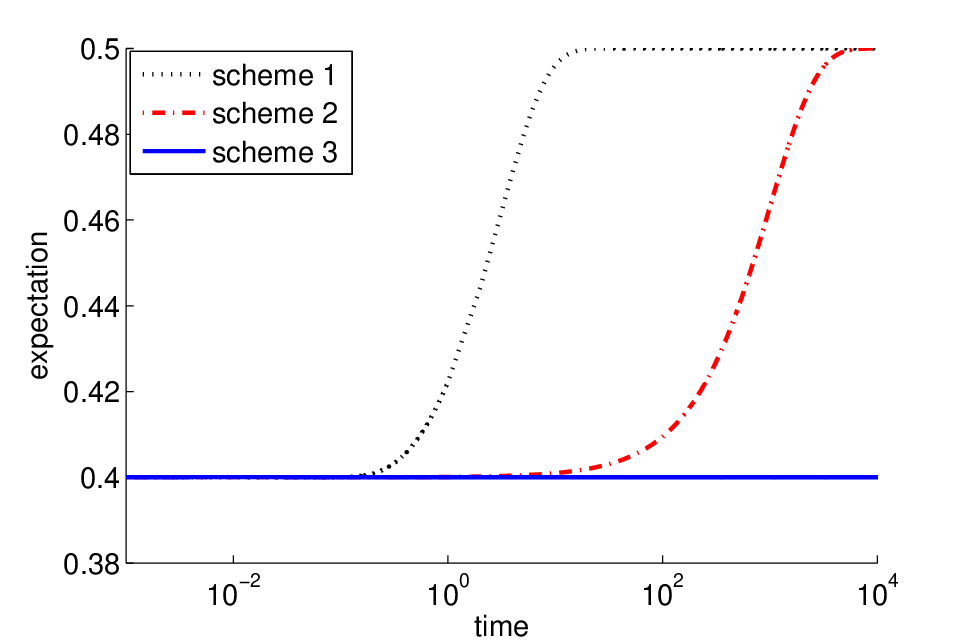}
\end{center}
\caption{Expectation produced by different schemes  with initial state $f(x,0) \sim
N(0.4,0.01^2)$ under the logarithm  time scale.}\label{exp}
\end{figure}

Fig. \ref{exp} shows that Scheme 3 preserves the expectation while Schemes 1 and 2 do not. This means only Scheme 3 can yield  the
\emph{complete solution}.
And it can be found that it takes a much longer time for Scheme 2 to achieve the steady state  than for Scheme 1.

 Table 1
presents the probability density and probability at the
boundaries ($f_{0}^{n}$,$f_{M}^{n}$ and $\frac{h}{2}f_{0}^{n}$,
$\frac{h}{2}f_{M}^{n}$), for Scheme 3 using different
space grid sizes ($h=\frac{1}{100}, \frac{1}{1000},
\frac{1}{10000}$) and different initial states ($f(x,0)\sim
N(0.4,0.01^2)$ and $N(0.7,0.01^2)$) at t=6 ($\tau=0.0001$). It
is shown that the probabilities at the boundary ends of the steady-state solution are independent of the space grid size $h$. This verifies that the Dirac delta singularities do develop at the boundary points. It is also shown that the probabilities at the boundary points are depended on the expectation of
the initial condition, and conservations of probability and expectation are always kept.    In Table 2, for fixed grid spacing $h=1/1000$ and two different initial states $f_0 \sim N(0.4, 0.01^2)$ and $f_0 \sim N(0.7,0.01^2)$, the time step $\tau$ is changed from $\frac 1 {10}$ to $\frac 1 {10000}$. The results shows that Scheme 3  can get the same steady state and keep the conservation of total property and expectation for different mesh ratio $\gamma=\frac{\tau}{h^2}$. This means that  the Scheme 3    for the
genetic drift problem  \eqref{con-dif0}, \eqref{eq03}, \eqref{bc1} and \eqref{ic} is stable and independent of the mesh ratio $\gamma$, i.e. unconditional stable.

%\begin{figure}[!htbp]
%\begin{center}
%\includegraphics[width=2.5in]{c1.eps}
%\includegraphics[width=2.5in]{c11.eps}\\
%\includegraphics[width=2.5in]{c2.eps}
%\includegraphics[width=2.5in]{c22.eps}\\
%\includegraphics[width=2.5in]{c3.eps}
%\includegraphics[width=2.5in]{c33.eps}\\
%\includegraphics[width=2.5in]{c4.eps}
%\includegraphics[width=2.5in]{c44.eps}
%\end{center}
%\caption{Numerical results of central finite volume method for genetic drift
%diffusion equation different initial states ($f(x,0) \sim
%N(0.7,0.01)$ (left column) and $f(x,0) \sim N(0.4,0.01)$ (right
%column) at different time steps (t=0.01, 0.05, 5, 6). The space
%and time step sizes are $h=0.0001,\ \tau =0.0001$. }
%\end{figure}

\begin{table}[!htbp]
\begin{center}
\begin{tabular}{|c|c|c|c|c|}
\hline
\multirow{1}{*}{space step
 $h$} & \multicolumn{4}{|c|}{$N(0.4,0.01^2)$} \\
\cline{2-5}
& $f_{0}^n$ & $f_{M}^{n}$ & $\frac{h}{2}f_{0}^{n}$ & $\frac{h}{2}f_{M}^{n}$  \\
\hline
$ 1/100$ &1.19999123e2 &7.99991237e1   & 0.59999562&  0.39999562\\
\hline
$1/1000$&1.19999115e3&   7.99991154e2&0.59999558 & 0.39999558\\
\hline
$1/10000$&1.19999115e4 & 7.99991146e3& 0.59999557&0.39999557\\
\hline
\multirow{1}{*}{space step
 $h$} & \multicolumn{4}{|c|}{$N(0.7,0.01^2)$} \\
\cline{2-5}
& $f_{0}^n$ & $f_{M}^{n}$ & $\frac{h}{2}f_{0}^{n}$ & $\frac{h}{2}f_{M}^{n}$  \\
\hline
$1/100$ & 5.99992332e1& 1.39999236e2 &0.29999617&  0.69999617\\
\hline
$1/1000$&  5.99992260e2&   1.39999226e3&0.29999613 &0.69999613\\
\hline
$1/10000$& 5.99992253e3 &  1.39999225e4&0.29999613& 0.69999613\\
\hline
\end{tabular}
\end{center}
\caption{Numerical results of Scheme 3 at the boundaries at steady state ($t=6$, time step $\tau=1/10000$)
with different initial states and space grid sizes.\label{tabel1} }
\end{table}

%\begin{table}[!htbp]
%\begin{center}
%\begin{tabular}{|c|c|c|c|c|c|c|c|c|}
%\hline
%\multirow{2}{*}{space step
% $h$} & \multicolumn{4}{|c|}{$N(0.4,0.1)$} & \multicolumn{4}{|c|}{$N(0.7,0.1)$} \\
%\cline{2-9}
%& $f_{0}^n$ & $f_{M}^{n}$ & $\frac{h}{2}f_{0}^{n}$ & $\frac{h}{2}f_{M}^{n}$ &$f_{0}^{n}$ & $f_{M}^{n}$ & $\frac{h}{2}f_{0}^{n}$ & $\frac{h}{2}f_{M}^{n}$  \\
%\hline
%$\frac{1}{100}$&  119.99 &79.999 & 0.5999& 0.3999 & 59.999 & 139.99 & 0.2999 & 0.6999\\
%\hline
%$\frac{1}{1000}$& 1199.9 &799.99 & 0.5999 & 0.3999 & 599.99 & 1399.9 & 0.2999 & 0.6999\\
%\hline
%$\frac{1}{10000}$&11999 & 7999.9 & 0.5999 & 0.3999 & 5999.9 & 13999 & 0.2999 & 0.6999\\
%\hline
%\end{tabular}
%\end{center}
%\caption{Numerical results of CFVM at the boundaries at steady state ($t=6$,$\tau=0.01$)
%with different initial states and space grid sizes.\label{tabel2} }
%\end{table}

\begin{table}[!htbp]
\begin{center}
\begin{tabular}{|c|c|c|c|c|}
\hline
\multirow{1}{*}{time step
 $\tau$} & \multicolumn{4}{|c|}{$N(0.4,0.01^2)$} \\
\cline{2-5}
& $f_{0}^n$ & $f_{M}^{n}$ & $\frac{h}{2}f_{0}^{n}$ & $\frac{h}{2}f_{M}^{n}$  \\
\hline
$1/10$  &1.19997448e3 &  7.99974480e2&0.59998724&  0.39998724\\
\hline
$1/100$ &1.19999005e3&   7.99990054e2&0.59999503&  0.39999502\\
\hline
$1/1000$&1.19999106e3&   7.99991058e2&0.59999553 & 0.39999553\\
\hline
$1/10000$&1.19999115e3&   7.99991154e2&0.59999558 & 0.39999558\\
 \hline
\multirow{1}{*}{time step
 $\tau$} & \multicolumn{4}{|c|}{$N(0.7,0.01^2)$} \\
\cline{2-5}
& $f_{0}^n$ & $f_{M}^{n}$ & $\frac{h}{2}f_{0}^{n}$ & $\frac{h}{2}f_{M}^{n}$  \\
\hline
$1/10$     &5.99977672e2 &1.39997767e3  &0.29998884&0.69998884\\
\hline
$1/100$    &5.99991298e2 &1.39999130e3  &0.29999565&0.69999565\\
\hline
$1/1000$   &5.99992177e2 &1.39999218e3  &0.29999609&0.29999608\\
\hline
$1/10000$  &5.99992260e2 &1.39999226e3  &0.29999613&0.69999613\\
\hline
\end{tabular}
\end{center}
\caption{Numerical results of Scheme 3 at the boundaries at steady state ($t=6$)
with same  space grid sizes $h= 1/1000$.\label{tabel2} }
\end{table}

%\begin{table}[!htbp]
%\begin{center}
%\begin{tabular}{|c|c|c|c|c|}
%\hline
%\multirow{1}{*}{time step
% $\tau$} &$f_{0}^n$ & $f_{M}^{n}$ & $\frac{h}{2}f_{0}^{n}$ & $\frac{h}{2}f_{M}^{n}$  \\
%\hline
%$\frac 1 {10}$& 1199.9787 &799.9787&  0.599989& 0.399989 \\
%\hline
%$\frac 1 {100}$&1199.9902 &799.9902 & 0.599995& 0.399995\\
%\hline
%$\frac 1 {1000}$&1199.9910& 799.9910 &0.599996 & 0.399996\\
%\hline
%$\frac 1 {10000}$&1199.9912 & 799.9912 &0.599996 & 0.399996\\
%\hline
%\end{tabular}
%\end{center}
%\caption{Numerical results of CFVM at the boundaries at steady state
%with same  space grid sizes $h=\frac 1 {1000}$ and initial state $f_0 \sim N(0.4,0.01^2)$.\label{tabel2} }
%\end{table}

\subsection{Analysis of Results}

\subsubsection{Analysis of   Schemes 1 and 2}   As shown in Figs.
\ref{up04}- \ref{s07}, in the results of the  first two schemes, the
values of steady-state solutions at boundaries $x=0$ and $x=1$ are
of the same height with different initial states. This is not
consistent with the steady state of the singular solutions given
in (\ref{stable}), and also does not satisfy the conservation of
the expectation. Generally, compared to schemes without upwind technique,
the upwind scheme is a better choice for the convection-diffusion
problem because it is more stable due to its intrinsic numerical
viscosity. It is not this case for the problem here. So we first check the effect of the numerical viscosity. To this purpose,   we consider a procedure of viscosity-vanishing. First, a infinitesimal diffusion is added to the corresponding steady-state  problem, then the limit behavior of the steady-state solution is investigated.
Consider the following problem, for a small  $\epsilon>0$,
\begin{equation}
\begin{cases}
-\frac{\mathrm {d}^2}{\mathrm
{d}x^2}\left((x(1-x)+\epsilon)f_\epsilon\right)=0, x\in (0,1),
\\
\frac{\mathrm {d}}{\mathrm
{d}x}\left((x(1-x)+\epsilon)f_\epsilon\right)\big{|}_{x=0,1}=0,
\end{cases}\label{eqv}
\end{equation}
with a constraint  $\int_0^1\ f_\epsilon\  \mathrm {d}x=1$.
Integrating the problem (\ref{eqv}),\ we get
\begin{equation}
f_\epsilon=\frac{a_\epsilon x+b_\epsilon}{x(1-x)+\epsilon}
\end{equation}
where $a_\epsilon$ and  $b_\epsilon$\ are related to $\epsilon $.\
From the boundary condition, \ we get $a_\epsilon=0$.\ The constraint on the total probability yields that
$$b_\epsilon=\frac {c_{\epsilon}^+}{\ln \frac{ c_{\epsilon}^++1/2}{c_{\epsilon}^+-1/2}},$$
where $c_{\epsilon}^+=\sqrt{\frac 1 4+\epsilon}.$
It is easy to verify that when $\epsilon \to 0$,
\begin{equation*}
f_\epsilon \to
\begin{cases}
0 & x\in (0,1) \\
\infty & x=0,1,
\end{cases}
\end{equation*}
(see Fig. 6 for $f_\epsilon$).

It is expected that  $f_\epsilon$ converges to
$\frac{1}{2}\delta(x)+\frac{1}{2}\delta(x-1)$ in the sense of distribution.\ Actually we have
the following theorem.

\begin{figure}[!htbp]
\begin{center}
\includegraphics[width=3.5in]{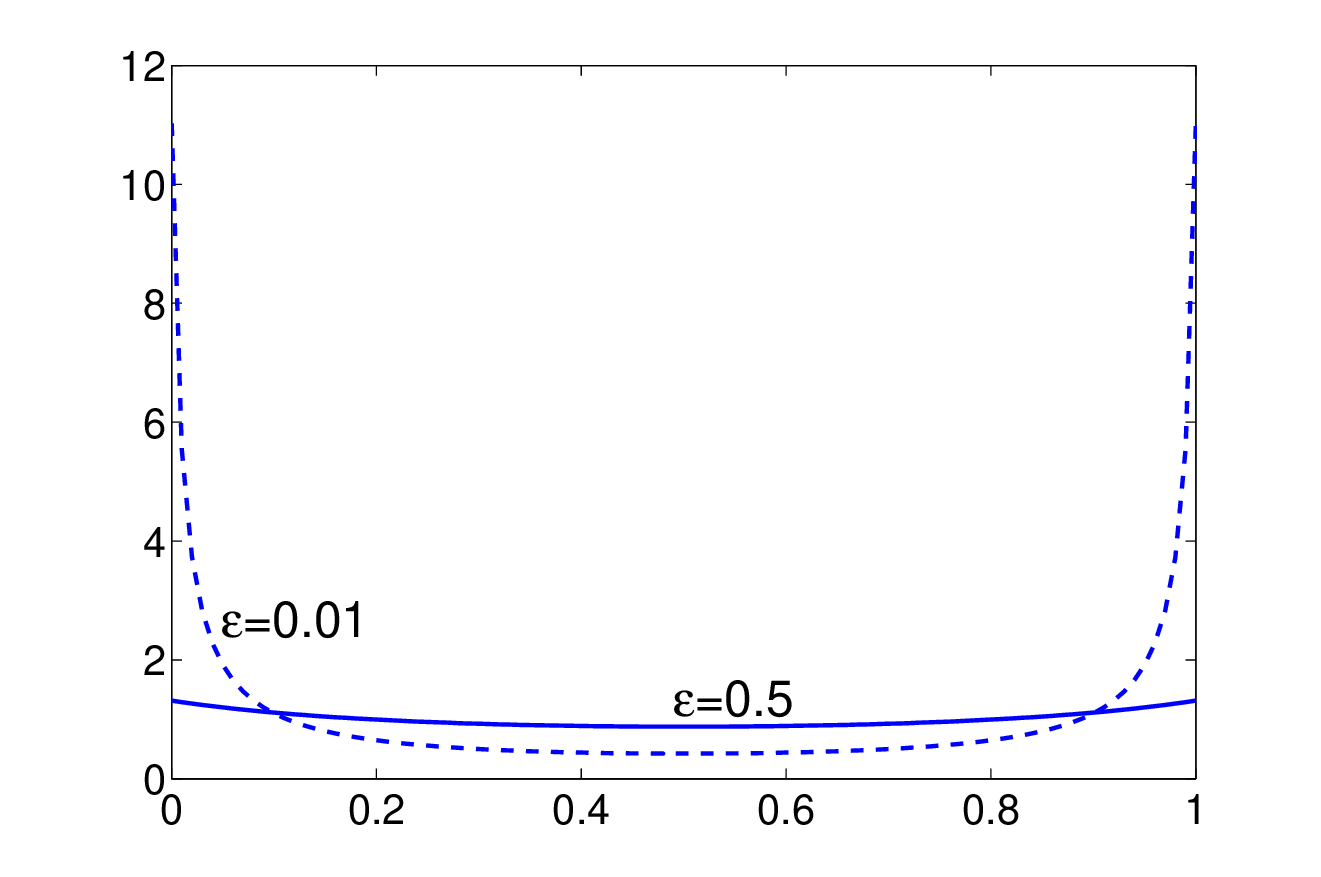}\label{fig4}
\caption{Two stable solution of genetic drift diffusion equation
added viscosity, $f_\epsilon$, with two different viscosities,
$\epsilon=0.5$, and $0.01$.}
\end{center}
\end{figure}

\begin{theorem} Take a zero extension of $\
f_\epsilon$ and still denote  by $\ f_\epsilon$,
\begin{equation*}
f_\epsilon =
\begin{cases}
\frac{b_\epsilon}{x(1-x)+\epsilon} & x \in [0,1]\\
0 & \mbox{otherwise},
\end{cases}
\end{equation*}
then
$  \underset {\epsilon \to 0} {\mathrm{lim} }
\int_{-\infty}^{+\infty}\ f_\epsilon \varphi\ \mathrm{d}x\ =\frac
{1}{2}\varphi(0)+\frac {1}{2}\varphi(1)$, \ $\ \forall \varphi \in \mathrm{C}_0^\infty (-\infty, \infty)$. \end{theorem}

 \noindent {\bf Proof.}\
\  For any $  \varphi \in \mathrm{C}_0^\infty (-\infty, \infty)$,\
$  \int_{-\infty}^{+\infty}\ f_\epsilon
\varphi\ \mathrm{d}x\ =\ \int_{0}^{1}\ f_\epsilon \varphi\
\mathrm{d}x $\\ First, choosing  $\varphi \in \mathrm{C}_0^\infty
(-\infty,\ \frac 1 2)$, we have
\begin{eqnarray*}
\underset {\epsilon \to 0} {\mathrm{lim} } <f_\epsilon,\varphi>\ &=&\ \underset {\epsilon \to 0} {\mathrm{lim} }  \int_{0}^{\frac 1 2}\ f_\epsilon \varphi\ \mathrm{d}x\\
&=&\ \underset {\epsilon \to 0} {\mathrm{lim} }  \int_{0}^{\frac 1 2}\  f_\epsilon\ (\varphi(x)-\varphi(0))+f_\epsilon \varphi(0)  \mathrm{d}x\\
&=&\ \underset {\epsilon \to 0} {\mathrm{lim} }  \int_{0}^{\frac 1 2}\  f_\epsilon\ \left(\varphi(x)-\varphi(0)\right)\  \mathrm{d}x+\varphi(0)\int_{0}^{\frac 1 2}f_\epsilon \  \mathrm{d}x\\
&=&\ \underset {\epsilon \to 0} {\mathrm{lim} }  \int_{0}^{\frac 1
2}\  f_\epsilon\ \left(\varphi(x)-\varphi(0)\right)\  \mathrm{d}x+\frac{1}{2}
\varphi(0).
\end{eqnarray*}
Now we   need to prove that  $\underset {\epsilon \to 0}
{\mathrm{lim} }  \int_{0}^{\frac 1 2} f_\epsilon
\left(\varphi(x)-\varphi(0)\right)\  \mathrm{d}x=0$. Denoting by $c_{\epsilon}^-=\sqrt{\frac 1 4-\epsilon}$ and $D(x)= x(1-x)$, we have
\begin{eqnarray*}
  \int_{0}^{\frac 1 2}f_\epsilon
\left(\varphi(x)-\varphi(0)\right)\  \mathrm{d}x&=& \int_{0}^{c_{\epsilon}^-}f_\epsilon\ \left(\varphi(x)-\varphi(0)\right)\
\mathrm{d}x\\&&+ \int_{c_{\epsilon}^-}^{\frac 1 2}
f_\epsilon\ \left(\varphi(x)-\varphi(0)\right)\
\mathrm{d}x\\
&\stackrel{\bigtriangleup}{=}&\mathrm{I}_1+\mathrm{I}_2,
\end{eqnarray*}
where
\begin{eqnarray*}
|\mathrm{I}_1 |&=&\left| b_\epsilon \int_{0}^{c_{\epsilon}^-}\frac{1}{D(x)+\epsilon}\ \left(\varphi(x)-\varphi(0)\right)\  \mathrm{d}x \right |\\
&\le& b_\epsilon M_1 \int_{0}^{c_{\epsilon}^-}\frac{x}{D(x)+\epsilon}\ \mathrm{d}x\ \ \ (M_1=\max\limits_{x\in [0,\frac{1}{2}]} |\varphi^\prime(x)|)\\
&\le&b_\epsilon M_1 c_{\epsilon}^-\frac{c_{\epsilon}^-}{c_{\epsilon}^-
(1-c_{\epsilon}^-)+\epsilon} \to 0, \qquad  \mbox{as}\ \epsilon \to 0
\end{eqnarray*}
and
\begin{eqnarray*}
|\mathrm{I}_2| &=&\left| b_\epsilon \int^{\frac 1 2}_{c_{\epsilon}^-}\ \frac{1}{D(x)+\epsilon}\ \left(\varphi(x)-\varphi(0)\right)\  \mathrm{d}x \right |\\
 &\le& M_2b_\epsilon\int^{\frac 1 2}_{c_{\epsilon}^-}\  \frac{1}{D(x)+\epsilon}\ \mathrm{d}x\  \  \ (
 M_2=\max\limits_{x\in[0,\frac{1}{2}]}|\varphi(x)-\varphi(0)|) \\
&=&b_\epsilon M_2 \frac{1}{2c_{\epsilon}^-}\ln \frac{c_{\epsilon}^+
-c_{\epsilon}^-+\frac 1 2}{c_{\epsilon}^++c_{\epsilon}^--\frac 1 2} \to 0, \ \ \mbox{as}\ \epsilon \to 0.
\end{eqnarray*}
Thus
\begin{equation}
\underset {\epsilon \to 0} {\mathrm{lim} } <f_\epsilon,\varphi>\
=\frac{1}{2}\varphi(0).
\end{equation}
Similarly, choosing  $\varphi \in \mathrm{C}_0^\infty (\frac 1 2,\
\infty)$,\  we get
\begin{equation}
\underset {\epsilon \to 0} {\mathrm{lim} } <f_\epsilon,\varphi>\
=\frac{1}{2}\varphi(1).
\end{equation}
Combining together,\  we have that  $\underset {\epsilon
\to 0} {\mathrm{lim} } <f_\epsilon,\varphi>\ =\frac
{1}{2}\varphi(0)+\frac{1}{2}\varphi(1)$, for $\forall \varphi(x)
\in \mathrm{C}_0^\infty(-\infty, \infty)$. The proof is completed. $\blacksquare$

 We  see that the long time
behavior of the original problem will be changed by any extra
infinitesimal viscosity  since the steady-state solution always is $\delta(x)/2+\delta(1-x)/2$, no matter what the initial state is. It  means that the conservation of expectation must be destroyed by any extra viscosity. This is consistent with the results showed in Figs. \ref{up04} and \ref{up07} for Scheme 1, since we know that the upwind scheme always introduces a first order viscosity into a central difference scheme.

However, numerical results in Figs. \ref{s04} and \ref{s07} show that the central difference scheme (Scheme 2) is stable for this convection-dominated problem and the
values of steady-state solutions at boundaries $x=0$ and $x=1$ are of the same height with different initial states, which also destroy the conservation of the expectation. Turning to Figs \ref{c04} and \ref{c07}, for the other central scheme (Scheme 3), we have stable and correct long time behavior. In the next subsection, we will prove the stability and long-time convergence for Scheme 3. Taking the difference between these schemes, we will get the answer for these different observations.

\begin{theorem} \label{diff_scheme}
Scheme 2: \eqref{eq8}, \eqref{eq6} and \eqref{scheme1} can be rewritten in the following form as
\be
\frac{f^n_i-f^{n-1}_i}{\tau}-\frac{D_{i+1}f^{n}_{i+1}-2D_if^n_i+D_{i-1}f^n_{i-1}}{h^2}+\Lambda=0
\ee
and Scheme 1: \eqref{eq8}, \eqref{eq6} and \eqref{scheme0} can be rewritten in the following form as
\be
\frac{f^n_i-f^{n-1}_i}{\tau}-\frac{D_{i+1}f^{n}_{i+1}-2D_if^n_i+D_{i-1}f^n_{i-1}}{h^2}+\Lambda+\tilde{\Lambda}=0,
\ee
with $$\Lambda=-\frac {h^2} {4}  \frac{f^n_{i+1}-2f^n_{i}+f^n_{i-1}}{h^2}, \ \ \tilde{\Lambda} =-\frac{h}{2}\frac{|b_{i+\frac 1 2}|f_{i+1}^n-(|b_{i+\frac 1 2}|+|b_{i-\frac 1 2})|f^n_{i}+|b_{i-\frac 1 2}|f^n_{i-1}}{h^2},$$
where $D(x) = x(1-x)$ and $b(x) = D_x = 1-2x$.
\end{theorem}

Theorem \ref{diff_scheme} tells us that Scheme 2 is just Scheme 3 plus a second order (of $O(h^2)$) viscosity term and for Scheme 1, another first order viscosity term is introduced in.  That is the reason why Scheme 2 is also stable and it takes a much longer time for Scheme 2 than for Scheme 1 to achieve the same wrong steady state.

{\bf Proof of Theorem \ref{diff_scheme}. }
  $\Lambda$ is the difference between Schemes 2 and 3. We have after direct computation that
 \be \ \ \ \
\Lambda &=& \frac 1 {h^2}\left((D_{i+1}-D_{i+\frac 1 2})f_{i+1}^n+(D_{i+\frac 1 2}+D_{i-\frac 1 2}-2D_i)f^n_i+(D_{i-1}-D_{i-\frac 1 2})f^n_{i-1}\right)\notag \\
        && - \frac{b_{i+\frac 1 2}f^n_{i+1}+(b_{i+\frac 1 2}-b_{i-\frac 1 2})f^n_{i}-b_{i-\frac 1 2}f^n_{i-1}}{2h}\notag \\
        &=&\frac{1}{h^2}\left((D_{i+1}-D_{i+\frac 1 2}-\frac h 2 b_{i+\frac 1 2} )f_{i+1}^n+(D_{i-1}-D_{i-\frac 1 2}+b_{i-\frac 1 2})f^n_{i-1}\right)\notag \\
        &&+\frac{1}{h^2}\left((D_{i+\frac 1 2}+D_{i-\frac 1 2}-2D_i-(b_{i+\frac 1 2}-b_{i-\frac 1 2}))f^n_i\right)\notag \\
        &=&\frac 1 {h^2}\left(\frac{h^2}{8}D_{xx} f^n_{i+1}+\frac{h^2}{8}D_{xx} f^n_{i-1}-\frac{h^2}{4}D_{xx} f^n_{i}\right)\notag \\
        &=&-\frac  { h^2}4 \frac{f^n_{i+1}-2f^n_{i}+f^n_{i-1}}{h^2},
\ee
where we have used the facts that $D_{xx}=-2$ and $\frac{d^{n}}{dx^n} D=0,~n>2$.

Note that $\tilde{\Lambda}$ is the difference between Schemes 1 and 2, which, as we have known, is a first order viscosity term introduced by a upwind scheme (\cite{Chen2015,Leveq}).  For completeness, we derive it again.  There are only 3 situations to check.

   \noindent \textbf{Case 1:} $0\le b_{i-\frac 1 2}< b_{i+\frac 1 2 }$.
\be
\tilde{\Lambda} &=& \frac{b_{i+\frac 1 2}f^n_{i+1}+(b_{i+\frac 1 2}-b_{i-\frac 1 2})f^n_i-b_{i-\frac 1 2}f^n_{i-1}}{2h}-\frac{ b_{i+\frac 1 2}f_{i+1}^n- b_{i-\frac 1 2}f_i^n}{ h}\nonumber\\
&=&\frac{-b_{i+\frac 1 2}f^n_{i+1}+(b_{i+\frac 1 2}-b_{i-\frac 1 2})f^n_i-b_{i-\frac 1 2}f^n_{i-1}}{2h}\nonumber\\
&=&-\frac{h}{2}\frac{b_{i+\frac 1 2}f^n_{i+1}-(b_{i+\frac 1 2}+b_{i-\frac 1 2})f^n_i+b_{i-\frac 1 2}f^n_{i-1}}{ h^2}.
\ee
\textbf{Case 2: } $  b_{i-\frac 1 2}< b_{i+\frac 1 2 }\le 0$.
\be
\tilde{\Lambda}
&=&-\frac{h}{2}\frac{-b_{i+\frac 1 2}f^n_{i+1}+(b_{i+\frac 1 2}+b_{i-\frac 1 2})f^n_i-b_{i-\frac 1 2}f^n_{i-1}}{ h^2}.
\ee
\textbf{Case 3: }$  b_{i-\frac 1 2}\le 0 \le  b_{i+\frac 1 2 }$.
\be
\tilde{\Lambda}
&=&-\frac{h}{2}\frac{b_{i+\frac 1 2}f^n_{i+1}-(b_{i+\frac 1 2}-b_{i-\frac 1 2})f^n_i-b_{i-\frac 1 2}f^n_{i-1}}{h^2}.
\ee
%
%\textbf{Case 4: }$b_{i+\frac 1 2}<0$ and $b_{i-\frac 1 2 }>0$,
%\be
%\tilde{\Lambda}
%&=&-\frac{h}{2}\frac{-b_{i+\frac 1 2}f^n_{i+1}-(-b_{i+\frac 1 2}+b_{i-\frac 1 2})f^n_i+b_{i-\frac 1 2}f^n_{i-1}}{2h}.
%\ee
Combining all together, we  get
$$ \tilde{\Lambda} =-\frac{h}{2}\frac{|b_{i+\frac 1 2}|f_{i+1}^n-(|b_{i+\frac 1 2}|+|b_{i-\frac 1 2})|f^n_{i}+|b_{i-\frac 1 2}|f^n_{i-1}}{h^2}.$$
The proof is completed. $\blacksquare$

\subsubsection{Analysis of   Scheme 3}
 In this section, we will prove the stability and long-time convergence for Scheme 3.
  \eqref{eq8}, \eqref{eq6} and \eqref{eq7} can be split into three independent parts.
For inner points, $i=1,\cdots, M-1,$
\begin{eqnarray}\label{inner}
\frac{f^{n+1}_i-f^{n}_i}{\tau}-\frac{D_{i+1}f^{n+1}_{i+1}-2D_{i}f^{n+1}_{i}+D_{i-1}f^{n+1}_{i-1}}{h^2}=0,~
\end{eqnarray}
with $D_i=x_i(1-x_i)$.
%In particular, due to $D_0=D_M=0$,
%\begin{eqnarray}
%\frac{f^{n+1}_1-f^{n}_1}{\tau}-\frac{D_{2}f^{n+1}_{2}-2D_{1}f^{n+1}_{1}}{h^2}=0, \\
%\frac{f^{n+1}_{M-1}-f^{n}_{M-1}}{\tau}-\frac{  a_{M-2}f^{n+1}_{M-2}-2D_{M-1}f^{n+1}_{M-1}}{h^2}=0,\label{boundary}
%\end{eqnarray}
For the boundary points, it yields,
\begin{eqnarray}
&& f^{n+1}_0=f^{n}_0+2D_1\gamma f^{n+1}_1,\\
&& f^{n+1}_M=f^{n}_M+2D_{M-1}\gamma f^{n+1}_{M-1},
\end{eqnarray}
with the mesh ratio $\gamma=\frac{\tau}{h^2}$.

Due to $D_0=D_M=0$, the unknowns at inner points $f^{n+1}_1,\cdots, f^{n+1}_{M-1}$  form a closed linear system which can be solved first.  Then  the boundary points $f^{n+1}_0, f^{n+1}_M$ can be updated by the inner points $f^{n+1}_1, ~f^{n+1}_{M-1}$ respectively.

In the following, we prove FVM keeps discrete total probability
and Scheme 3 also preserves the expectation. The discrete total
probability and expectation at step $n$ are defined as follows.
\begin{eqnarray}
P_n&=&\frac{h}{2}f_{0}^{n}+\frac{h}{2}f_{M}^{n}+\sum_{i=1}^{M-1}f_{i}^{n}h,\\
E_n&=&\frac{h}{2}x_0f_{0}^{n}+\frac{h}{2}x_{M}f_{M}^{n}+\sum_{i=1}^{M-1}x_if_{i}^{n}h,
\end{eqnarray}
\begin{lemma}\label{conserve} For the
genetic drift problem  \eqref{con-dif0}, \eqref{eq03}, \eqref{bc1}  and \eqref{ic}, the finite volume method
 \eqref{eq8}, \eqref{eq6}  keep the discrete total
probability, $P_{n+1}=P_{n}$. Furthermore, the central finite volume method \eqref{eq7}
 yields a complete solution, which preserves the discrete expectation,  $E_{n+1}=E_n$.
\end{lemma}

\noindent {\bf Proof.} By the definition of $P_n$, it yields,
\begin{equation}
P_{n+1}-P_{n}=\frac{h}{2}(f_{0}^{n+1}-f_{0}^{n})+\frac{h}{2}(f_{M}^{n+1}-f_{M}^{n})+\sum_{i=1}^{M-1}(f_{i}^{n+1}-f_{i}^{n})h.
\end{equation}
Using \eqref{eq8} and \eqref{eq6}, we have
\begin{eqnarray}P_{n+1}-P_{n}&=&-\tau
\sum_{i=1}^{M-1}(j_{i+\frac{1}{2}}^{n+1}-j_{i-\frac{1}{2}}^{n+1}) -\tau j_{\frac{1}{2}}^{n+1}+\tau
j_{M-\frac{1}{2}}^{n+1}=0.\end{eqnarray}
On the other hand,
\begin{equation}
E_{n+1}-E_{n}=\frac{h}{2}x_0(f_{0}^{n+1}-f_{0}^{n})+\frac{h}{2}x_M(f_{M}^{n+1}-f_{M}^{n})+\sum_{i=1}^{M-1}x_i(f_{i}^{n+1}-f_{i}^{n})h.
\end{equation}
Using \eqref{eq8} and \eqref{eq6}, we have
\begin{eqnarray}E_{n+1}-E_{n}&=&-\tau
\sum_{i=1}^{M-1}x_i(j_{i+\frac{1}{2}}^{n+1}-j_{i-\frac{1}{2}}^{n+1}) +\tau x_M j_{M-\frac{1}{2}}^{n+1}\ \ \ \notag \\
&=&h\tau\sum_{i=1}^{M-1}j_{i+\frac{1}{2}}^{n+1}.\end{eqnarray}
Combining with \eqref{eq7}, we get
\begin{equation}
\sum_{i=1}^{M-1}j_{i+\frac{1}{2}}^{n+1}=0.
\end{equation}
Thus, for Scheme 3, $E_{n+1}=E_{n}$. The proof is completed.  $\blacksquare$

\begin{lemma}\label{l2decay}
Let $f_i^n$ be the  solution of   Scheme 3: \eqref{eq8}, \eqref{eq6} and \eqref{eq7}. Then,
\begin{equation}
\sum_{i=1}^{M-1} x_if_i^{n+1}h\leq \sum_{i=1}^{M-1}f_i^{n+1}h\leq\frac {e^{-\frac{\pi^2t_{n+1}}{4}}}{4h(1-h)}\sum_{i=1}^{M-1}f_i^{0}h,~ t_{n+1}=(n+1)\tau.\nonumber
\end{equation}
\end{lemma}
\noindent {\bf Proof.}
Since $D_i=x_i(1-x_i)> 0$ for $i=1,\cdots,M-1$, by setting $v_i^n=D_if_i^n$,
 Scheme 3   for inner points \eqref{inner}  can be rewritten as
\begin{equation}\label{ftov}
\frac 1 {D_i}\frac{v_i^{n+1}-v_i^n}{\tau}-\frac{v_{i+1}^{n+1}-2v_i^{n+1}+v_{i-1}^{n+1}}{h^2}=0.
\end{equation}
It is obviously that the discrete maximum principle is valid. So we have
\begin{equation} \label{positive}
v_i^n \ge 0, \ \ f_i^n \ge 0, \ \ i=1, \cdots, M-1, \ n>1,
\end{equation}
since the initial value is nonnegative.

Multiplying by $v_i^{n+1}h$ on both sides of \eqref{ftov}, summing from $i=1$ to $M-1$, using H\"{o}lder inequality and $D(x)\leq \frac 1 4$, it yields that
\begin{equation}\label{est1}
2\frac {\|v_i^{n+1}\|^2_h}{\tau}-\sum_{i=1}^{M-1}\frac{v_{i+1}^{n+1}-2v_i^{n+1}+v_{i-1}^{n+1}}{h^2}v_i^{n+1}h\leq 2\frac{\|v_i^{n}\|^2_h}{\tau},
\end{equation}
where $\displaystyle\|v^n\|_{h} =\left(\sum_{i=1}^{M-1}|v_i^n|^2h\right)^{1/2}$ is the discrete $L^2$ norm of $v_n$.

Denote by ${\displaystyle \Delta_h w_i = \frac{w_{i-1}-2w_i+w_{i+1}}{h^2}}$ the discrete Laplacian. We know that the minimum eigenvalue of the problem $$-\Delta_h w_i  = \lambda w_i, \ \ i=1, \cdots, M-1, \ \ w_0 = w_M = 0,$$
is $\lambda_0=\frac 4 {h^2}\sin^2(\frac{\pi}{2M})$. Then from \eqref{est1},
\begin{equation}
2\frac {\|v_i^{n+1}\|^2_h}{\tau}+\lambda_0 \|v^{n+1}\|_h^2\leq 2\frac{\|v_i^{n}\|^2_h}{\tau}.
\end{equation}
%\begin{equation}A=\frac 1 {h^2}\left[\begin{array}{ccccccc}
%2 & -1 &   & & & & \\
%-1& 2 &-1  &  & & & \\
% &\ddots&\ddots&\ddots& & & \\
% &       &     &       &-1&2&-1\\
%  &       &     &       & &-1&2
%  \end{array}\right].\end{equation}
So we  get the estimate
\begin{equation}\label{vdecay}
\|v^{n+1}\|_{h}^2\leq\left( \frac 1 {1+\frac {\tau \lambda_0}{2}}\right) \|v^n\|_h^2\leq e^{-\frac{\pi^2t_{n+1}}{2}}\|v^0\|_h^2,~ t_{n+1}=(n+1)\tau.
\end{equation}
By H\"{o}lder inequality, above formula means,
\begin{equation}
\sum_{i=1}^{M-1} v_i^{n+1}h\leq (\sum_{i=1}^{M-1}|v_i^{n+1}|^2h)^{1/2}=\|v^{n+1}\|_h\leq e^{-\frac{\pi^2t_{n+1}}{4}}\|v^0\|_h.
\end{equation}
For fixed grid spacing $h$,  $h(1-h)\leq D_i\leq \frac 1 4$,  $i=1,\cdots, M-1$.
Thanks to $v_i^{n+1}=D_if_i^{n+1}$,  $f_i^{n+1}$  can be estimated as follows:
\begin{equation}\label{massdecay}
\sum_{i=1}^{M-1}f_i^{n+1}h\leq\frac {e^{-\frac{\pi^2t_{n+1}}{4}}}{4h(1-h)}\sum_{i=1}^{M-1}f_i^{0}h.
\end{equation}
At the same time, since $x\in (0,1)$, by \eqref{massdecay}, we can directly obtain,

\begin{equation}\label{expdecay}
\sum_{i=1}^{M-1} x_if_i^{n+1}h\leq \sum_{i=1}^{M-1}f_i^{n+1}h\leq\frac {e^{-\frac{\pi^2t_{n+1}}{4}}}{4h(1-h)}\sum_{i=1}^{M-1}f_i^{0}h.
\end{equation}
The proof is completed.  $\blacksquare$

\begin{theorem}   Scheme 3: \eqref{eq8}, \eqref{eq6} and \eqref{eq7} is unconditionally stable and for fixed $h$ and $\tau$ and   we have, as $n\rightarrow \infty$, that
 \begin{equation}
\left\{\begin{array}{l}f_i^n\rightarrow 0,~\mbox{for}~i=1,\cdots,M-1,\\
\frac h 2 f_0^{n}+\frac h 2 f_M^{n}\rightarrow P_0,\\
 \frac h 2 x_Mf_M^{n} \rightarrow E_0,\end{array}\right.
 \end{equation}
 i.e.,
 \begin{equation}
\left\{\begin{array}{l}
\frac h 2 f_0^{n} \rightarrow   (P_0- E_0 ) \approx 1-p,\\
\frac h 2 f_M^{n} \rightarrow    E_0 \approx p.\end{array}\right.
 \end{equation}

\end{theorem}
\noindent\textbf{Proof.} The unconditional stability is ensured by the discrete  maximum principle. The long-time convergence comes from   Lemmas \ref{conserve}, \ref{l2decay} and nonnegative property of the solution \eqref{positive}. The proof is completed. $\blacksquare$

\section{Conclusion and discussion}
\setcounter{equation}{0}
 We have considered three different numerical schemes for genetic drift diffusion equation. Schemes 1 and 2 discritize the flux in convection and diffusion form using upwind and central difference methods, respectively. Scheme 3 discritizes the flux as a whole using cental difference method. Numerical experiments show that all the schemes are stable and the first two schemes give an artificial steady state solution while   Scheme 3  presents the correct one.
Analysis shows that Scheme 3 preserves both of total probability and expectation, so it yields the complete solution. We prove that any extra infinitesimal diffusion leads to a same  steady state. We find that Scheme 2, also a central scheme, introduces a second order   numerical viscosity term to Scheme 3, while Scheme 1, as a upwind one, introduces another first order numerical viscosity term to Scheme 2. Therefore, both  Schemes 1 and 2 yield the same artificial long time behaviors of the numerical solutions, but it takes different time for them to achieve the steady state since they introduce different scales of viscosity terms. Some interesting observations are that a central scheme could be unconditional stable for a convection-dominated problem and a central scheme could also introduce numerical viscosity, which are beyond the common understanding of the convection-diffusion community.  All the complexity  comes from the diffusion degeneration and convection dominating. For this kind of problem, the numerical methods must be carefully chosen. Any methods with intrinsic numerical viscosity must be  avoided. This means that most of stable methods for convection-diffusion problem don't work since they achieve the stability by introducing the numerical viscosity (\cite{Leveq}).

Discontinuous Galerkin (DG) finite element \cite{Cockburn2000}  method  and a number of variants, such as  local discontinuous Galerkin (LDG)\cite{Cockburn2007, Cockburn1998},   Interior Penalty discontinuous Galerkin (IPDG)\cite{Arnold1982},  central discontinuous Galerkin (CDG) \cite{Liu2007}, are  also the  main tools to dealing with convection-diffusion equation.  We can get an equivalent form of scheme 3 by subtly designing the numerical flux of LDG and get the right results. The numerical flux should keep the inner nodes independent on the boundary nodes. Otherwise, DG methods also may not yield  the right results.

If a population or species of organisms typically includes multiple alleles at each locus among various individuals, which is called multiple alleles, the problem will be a high dimension problem (\cite{Kimura1955,Tran2012,Waxman2009}). It is a great challenge to find a  \emph{complete solution} since the singularity will always be developed on the boundary surface rather than only two points for 1-D case. The numerical method for the multiple alleles include the fixation phenomena will be  our future work.
\\
\\
\noindent
{\bf Acknowledgments}
 The authors benefitted a great deal from
discussions with   David Waxman. This work  is supported
in part by NSF of China under the grants 11271281, 11301368  and 91230106. Minxin Chen, Chun Liu and Shixin Xu are partially supported by NSF grants DMS- 1159937, DMS-1216938 and DMS-1109107.

\end{document}